\documentclass[reqno,12pt]{amsart}
\usepackage{amsmath, amssymb, amsthm,amsfonts} 
\usepackage{bbm,bm}
\usepackage{graphicx}
\usepackage{url}
\usepackage{epstopdf}
\usepackage[a4paper,bindingoffset=0.2cm,left=1cm,right=1cm,top=2.5cm,bottom=2cm,footskip=.8cm]{geometry}
\usepackage{subcaption}
\usepackage{textcomp}
\usepackage[algo2e,ruled,ruled,vlined,linesnumbered]{algorithm2e}
\usepackage{multicol}
\usepackage{multirow}
\usepackage{wrapfig}
\usepackage{diagbox}
\usepackage{floatrow}
\newfloatcommand{capbtabbox}{table}[][\FBwidth]
\usepackage[dvipsnames]{xcolor}
\usepackage{subfiles}
\usepackage{import}
\usepackage[T2A]{fontenc}
\usepackage[utf8]{inputenc}
\usepackage{booktabs}
\usepackage{hyperref}
\usepackage[english]{babel}

\usepackage{amsbsy,enumerate}
\usepackage{graphicx}
\usepackage{comment}
\usepackage{mathrsfs} 
\usepackage{xcolor}
\usepackage{mathtools}

\newcommand{\bs}{\boldsymbol}

\newcommand{\vb}{\vspace{5mm}}
\renewcommand{\hat}{\widehat}

\DeclareMathOperator*{\argmin}{arg\,min}

\allowdisplaybreaks

\makeatletter
\newtheorem*{rep@theorem}{\rep@title}
\newcommand{\newreptheorem}[2]{%
\newenvironment{rep#1}[1]{%
 \def\rep@title{#2 \ref{##1}}%
 \begin{rep@theorem}}%
 {\end{rep@theorem}}}
\makeatother

\makeatletter
\renewcommand*\env@matrix[1][\arraystretch]{%
  \edef\arraystretch{#1}%
  \hskip -\arraycolsep
  \let\@ifnextchar\new@ifnextchar
  \array{*\c@MaxMatrixCols c}}
\makeatother

\newtheorem{theorem}{Theorem}
\newreptheorem{theorem}{Theorem}
\newtheorem{lemma}{Lemma}

\newtheorem{remark}{Remark}

\newtheorem{proposition}{Proposition}
\newreptheorem{proposition}{Proposition}

\newcommand{\labda}{\lambda} 

\newcommand{\omegaT}{\omega^{\mathrm{T}}}
\newcommand{\omegaI}{\omega^{\mathrm{I}}}
\newcommand{\omegaW}{\omega^{\mathrm{W}}}
\newcommand{\scv}{\textsc{scv}}

\newcommand{\Tour}{\boldsymbol{\sigma}} 
\newcommand{\tour}{\sigma} 
\newcommand{\IA}{{\bs x}}
\newcommand{\EX}{\mathbb{E}\,} 

\newcommand{\objfunnew}{\mathscr{L}(\Tour, \IA)}
\newcommand{\objfunhyb}{\mathscr{L}_{\text{hyb}}}

\newcommand{\Bharti}[1]{\textcolor{OliveGreen}{[Bharti: #1]}}

\begin{document}

\title[Integrated routing and appointment scheduling]
{A queueing-based approach for \\integrated routing and appointment scheduling}

\author{Ren\'e Bekker, Bharti Bharti, Leon Lan, Michel Mandjes}

\maketitle

\begin{abstract}
    This paper aims to address the integrated routing and appointment scheduling (RAS) problem for a single service provider.
    The RAS problem is an operational challenge faced by operators that provide services requiring  home attendance, such as grocery delivery, home healthcare, or maintenance services.
    While considering the inherently random nature of service and travel times, the goal is to minimize a weighted sum of the operator's travel times and idle time, and the client's waiting times.
    To handle the complex search space of routing and appointment scheduling decisions, we propose a queueing-based approach to effectively deal with the appointment scheduling decisions.
    We use two well-known approximations from queueing theory: first, we use an approach based on phase-type distributions to accurately approximate the objective function, and second, we use an heavy-traffic approximation to derive an efficient procedure to obtain good appointment schedules.
    Combining these two approaches results in a fast and sufficiently accurate hybrid approximation, thus essentially reducing RAS to a routing problem.
    Moreover, we propose the use a simple yet effective large neighborhood search metaheuristic to explore the space of routing decisions.
    The effectiveness of our proposed methodology is tested on benchmark instances with up to 40 clients, demonstrating an efficient and accurate methodology for integrated routing and appointment scheduling.
    

\vb

\noindent
{\sc Affiliations.} \emph{Bharti Bharti} and \emph{Michel Mandjes} are with the Korteweg-de Vries Institute for Mathematics, University of Amsterdam, Science Park 904, 1098 XH Amsterdam, the Netherlands. MM is also with Mathematical Institute, Leiden University, The Netherlands, and Amsterdam Business School, Faculty of Economics and Business, University of Amsterdam, Amsterdam, the Netherlands. This research was supported by the European Union’s Horizon 2020 research and innovation programme under the Marie Skłodowska-Curie grant agreement no.\ 945045, and by the NWO Gravitation project NETWORKS under grant no.\ 024.002.003. \includegraphics[height=1em]{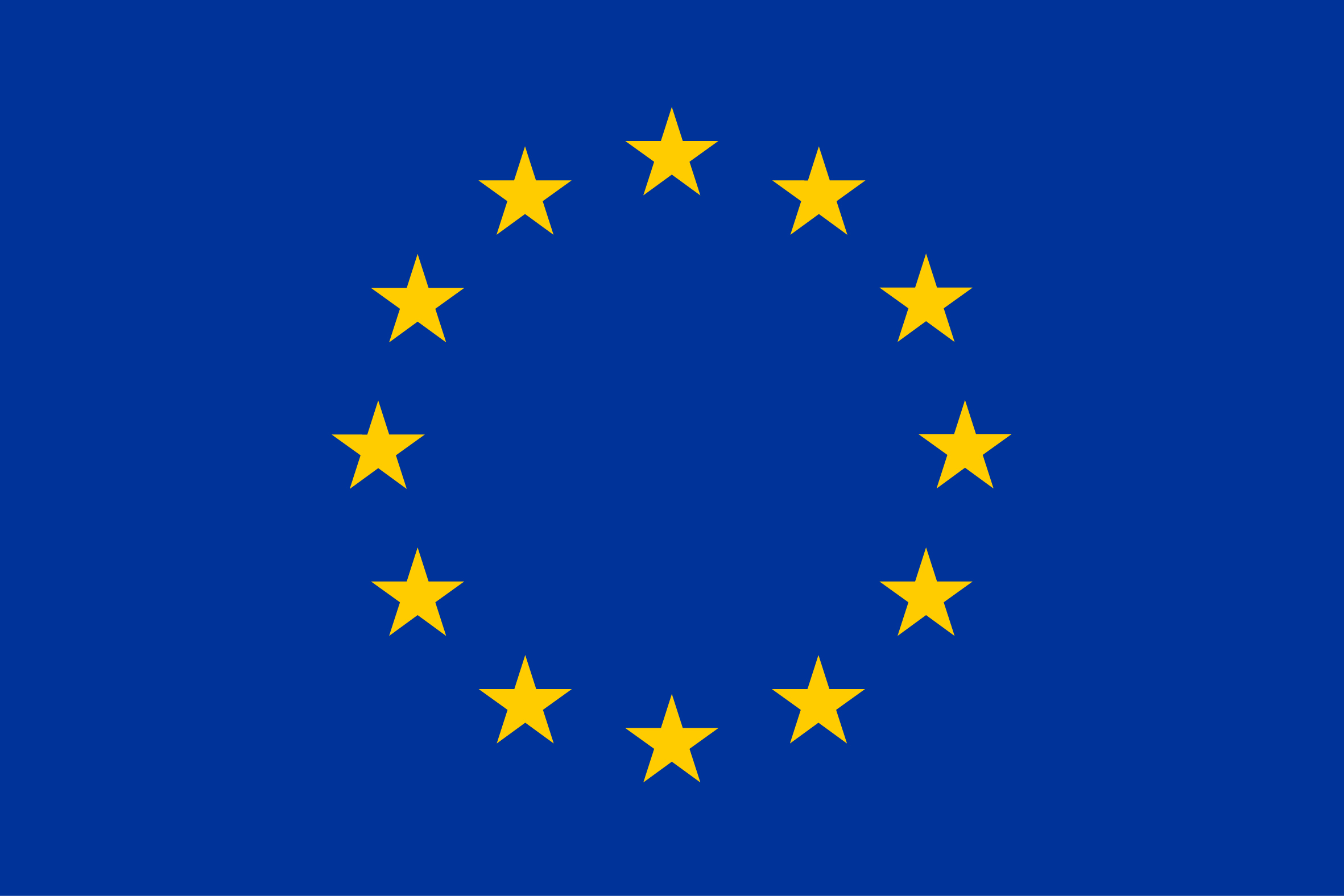}

\noindent
\emph{Ren\'e Bekker} and \emph{Leon Lan} are with the
Department of Mathematics,
Vrije Universiteit,
De Boelelaan 1111,
1081 HV Amsterdam,
the Netherlands. 

\vb

\noindent
\end{abstract}

\newpage

\section{Introduction}

The market for providing services at home is growing rapidly and encompasses a wide range of services, such as food delivery, home healthcare and wellness, maintenance and appliance repair, household activities, and more. A specific challenge for home attendance services is the combination of an efficient home service process combined with good quality of service (e.g. measured as `right on time' home attendance).
Hence, the operational challenge that the service provider is facing lies in determining, for a given set of clients with known locations, the corresponding route and appointment times.
This combined routing and appointment scheduling problem has become a critical issue, as inefficient scheduling and sub-optimal resource allocation can have significant cost implications. 
Studies \cite{srinivas2018} indicate that these inefficiencies cost the US healthcare system more than $\$150$ billion annually, resulting in reduced provider productivity, delayed access to care, and increased burden on the healthcare system.
Similarly, in the transportation sector, improper scheduling and routing have been associated with increased fuel consumption, carbon emissions, and total logistics delivery costs \cite{li2020}. 
The significance of this issue extends beyond financial aspects, encompassing enhanced service quality, improved client satisfaction, and the increasing pressure and competition faced by delivery companies to ensure on-time performance \cite{camacho2006}.
Therefore, it is important to study this problem, so as to develop effective scheduling strategies and solutions.

\begin{figure}
    \centering
 \begin{subfigure}[b]{0.75\textwidth}
   \includegraphics[width=1\linewidth]{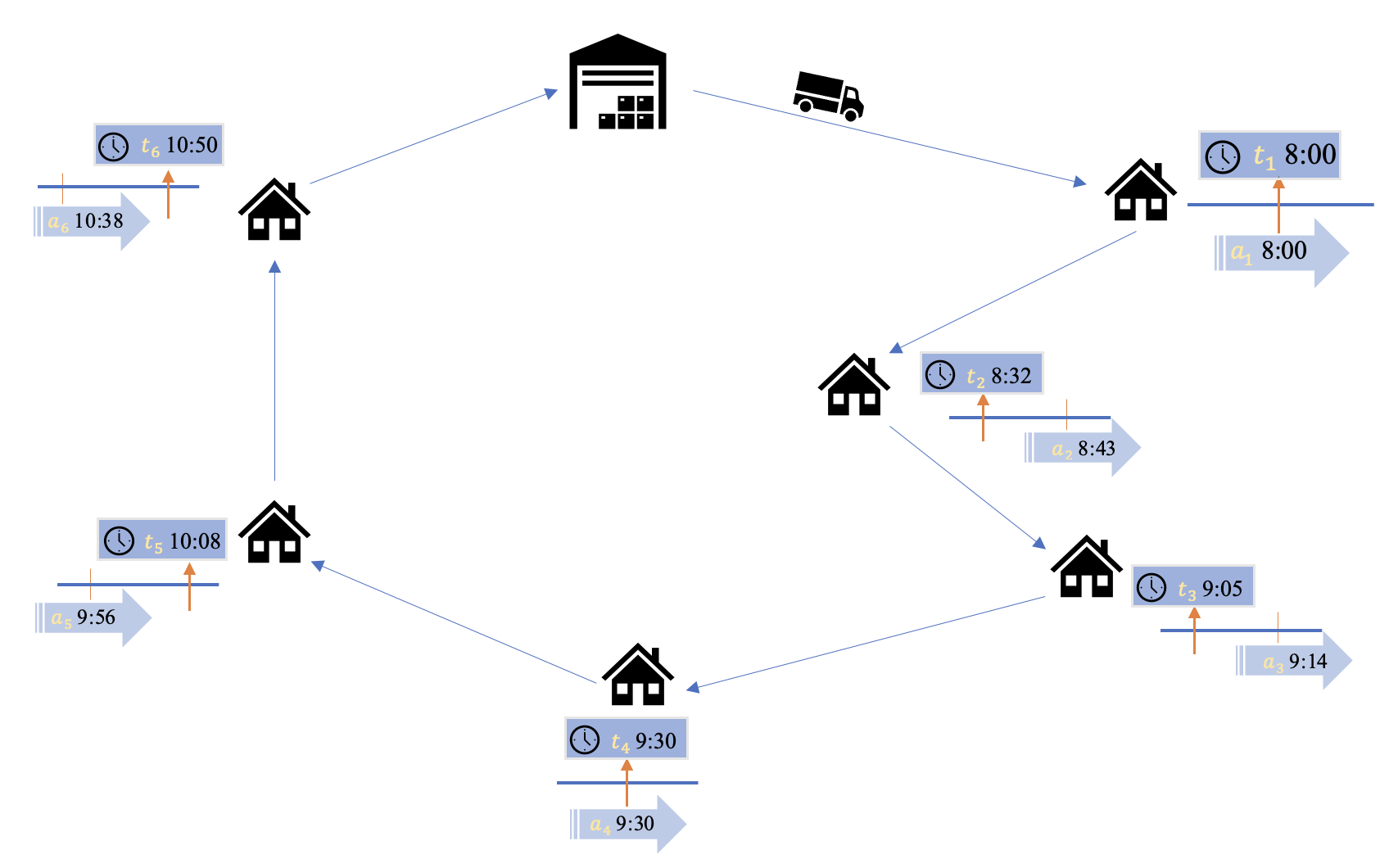}
   \caption{}
   \label{fig:Ng1} 
\end{subfigure}

 \begin{subfigure}[b]{0.75\textwidth}
   \includegraphics[width=1\linewidth]{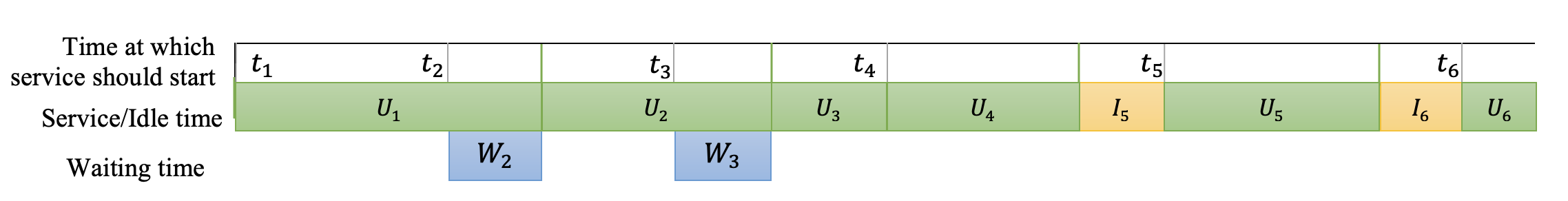}
   \caption{}
\end{subfigure}
\caption{
\label{F1} (A) shows a service provider with a route consisting of six locations. Here, $t_j$ and $a_j$ are the appointment time and arrival time at the $j^{\rm th}$ location, respectively. (B) shows the relation between `extended service time' $U_j$, idle time $I_j$ of the service provider before service of $j^{\rm th}$ client can start, and waiting time $W_j$ of the $j^{\rm th}$ client.
}
\end{figure}

The main objective of this paper is to develop a unified approach to solve the integrated routing and appointment scheduling (RAS) problem; an illustration of the RAS problem is depicted in Figure~\ref{F1}.
As opposed to what is usually considered in the literature, the setup we consider is intrinsically random.
More precisely, to do justice to their inherently random nature, we model the service times (the times spent at the clients' locations) and travel times (the time it takes to travel between two subsequent locations) as random variables.
The service provider is required to visit a predetermined set of locations before returning to their initial location, e.g., a warehouse.
The clients residing at these locations are provided with appointment times.

Importantly, the service provider and the clients have conflicting interests.
The service provider wishes to minimize their travel time and, in addition, they want to avoid arriving at the client's location before the corresponding appointment time resulting in idle time of the service provider. Instead, the clients prefer short waiting times. 
Therefore, in our approach, we aim to minimize these three elements so as to strike a proper balance between the interests of both the service provider and the clients.
Even in a setting without any randomness, this problem is intrinsically hard: the routing decision leads to the traveling salesman problem (TSP), which in itself is an NP-hard problem. 
A major challenge lies in the randomness that we wish to incorporate: it is not evident how to accurately and efficiently evaluate the objective function for given distributions of travel times and service times.


\vspace{3mm}

{\it Problem structure, benchmark algorithms.} 
Routing and appointment scheduling are mature subfields of operations research, in each of which a wealth of well-performing algorithms is available. 
In appointment scheduling, an excellent heuristic for the ordering of clients is the so-called {\it smallest variance first} (SVF) rule, which sequences clients in order of increasing variance of their service durations.
While this rule is technically not optimal \cite{ahmadi2017}, it has been formally proven that in specific asymptotic regimes it {\it is} \cite{cook2015}. 
The informal backing for SVF's near-optimality lies in the fact that one wishes to introduce as little as possible uncertainty into the system. 
If one schedules relatively deterministic jobs first, they hardly do any harm to their successors.

\begin{figure}[!ht]
  \centering
  \begin{subfigure}[b]{0.32\textwidth}
    \includegraphics[width=\textwidth]{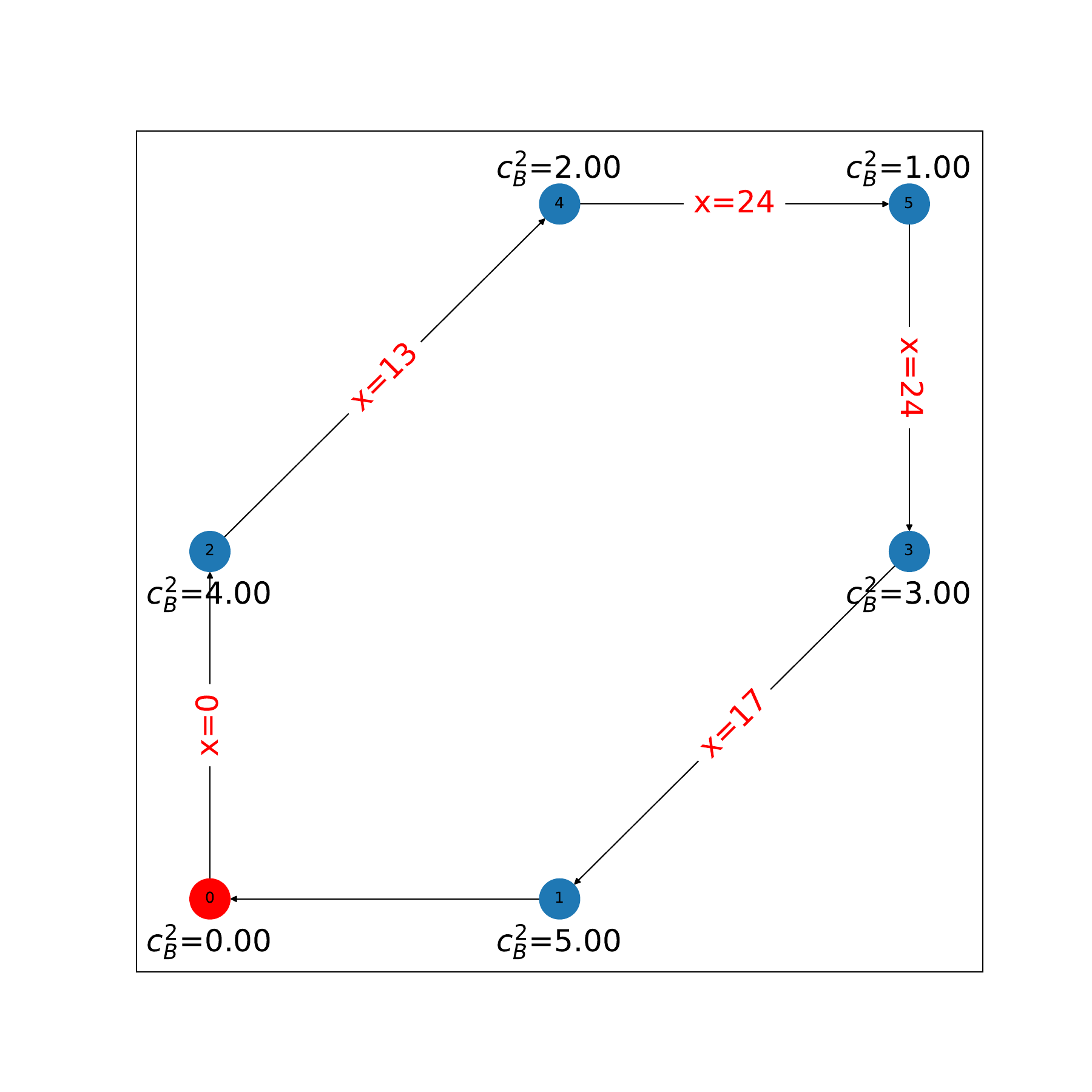}
    \caption{TSP}
    \label{fig:subfig1}
  \end{subfigure}
  \begin{subfigure}[b]{0.32\textwidth}
    \includegraphics[width=\textwidth]{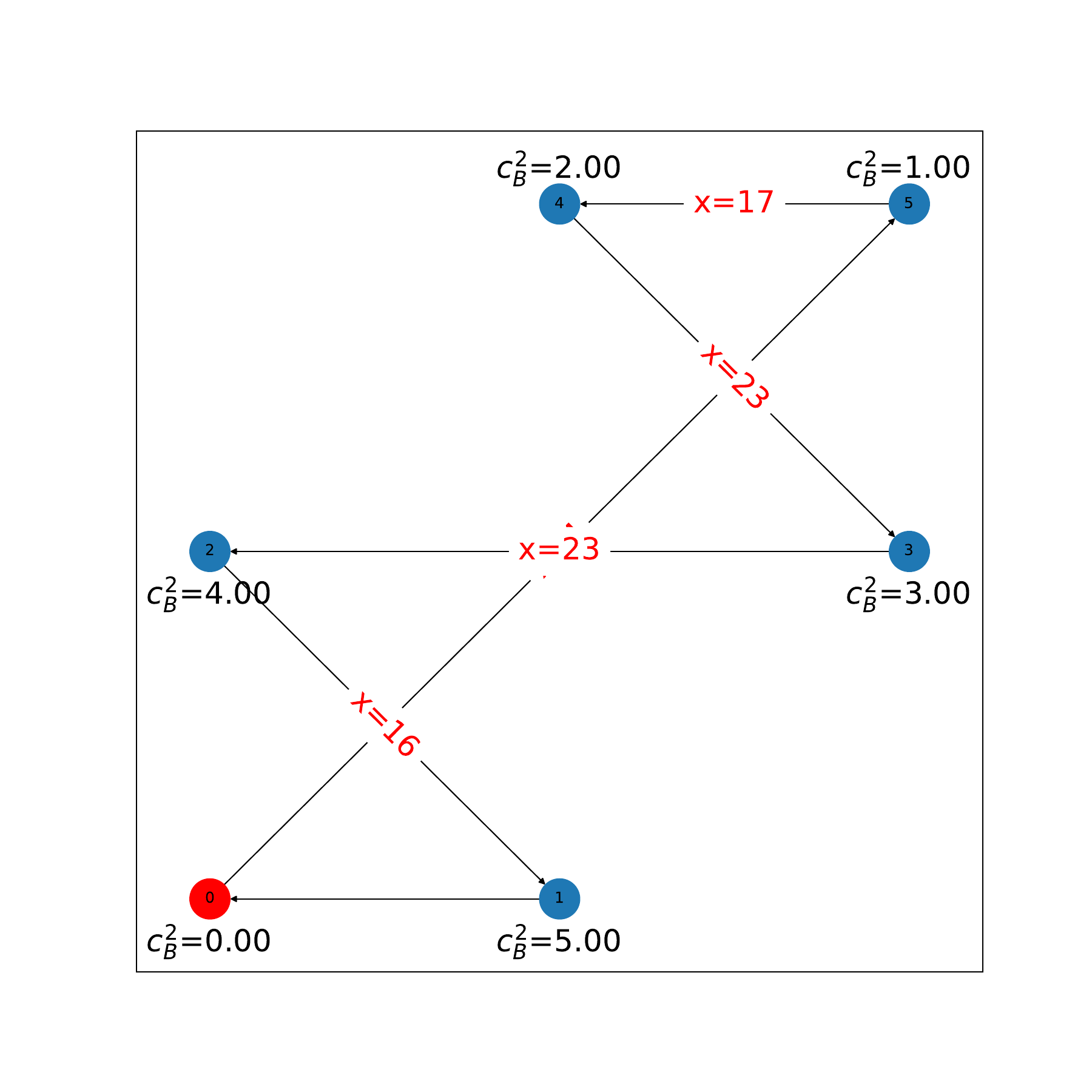}
    \caption{SVF}
    \label{fig:subfig2}
  \end{subfigure}
  \begin{subfigure}[b]{0.32\textwidth}
    \includegraphics[width=\textwidth]{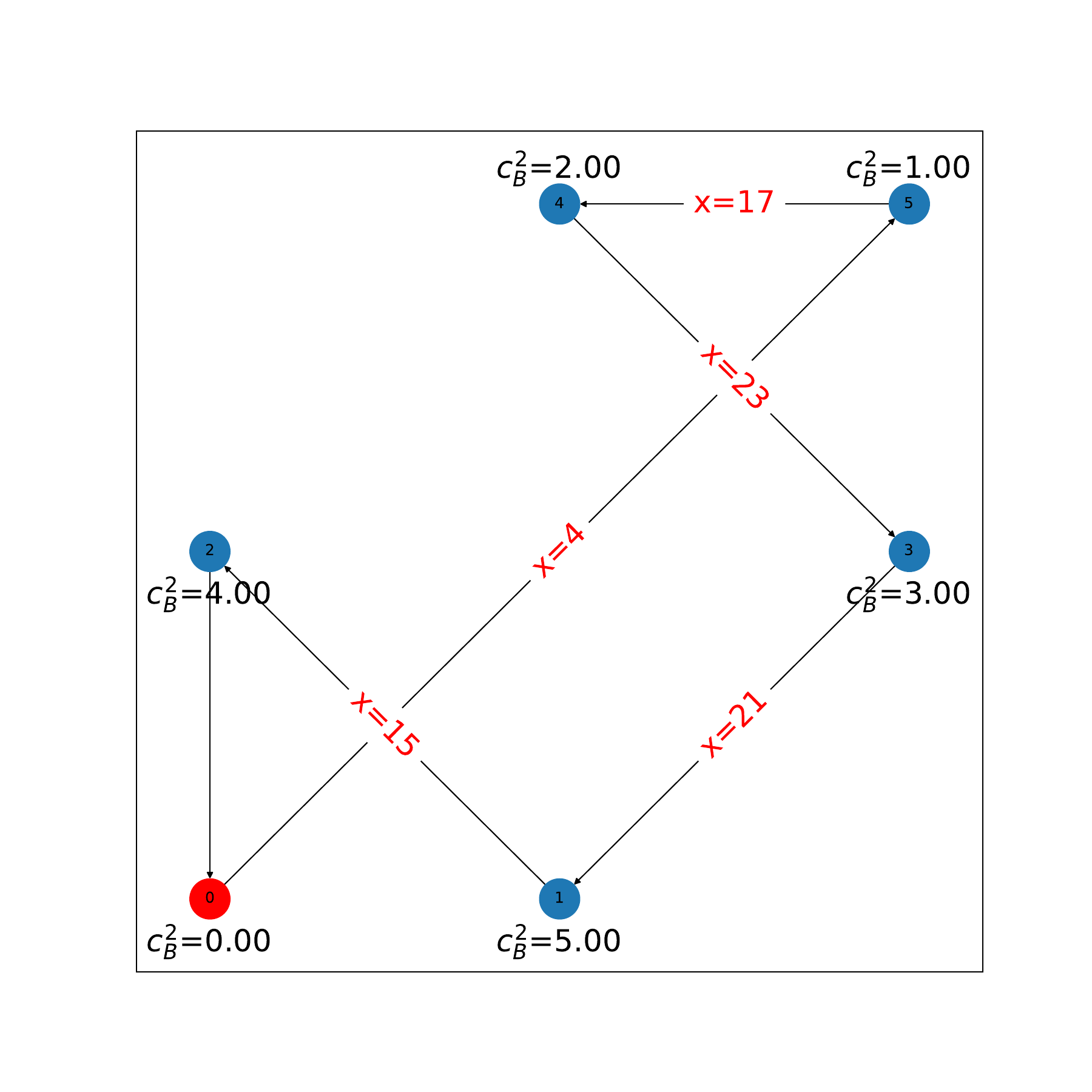}
    \caption{Optimal}
    \label{fig:subfig3}
  \end{subfigure}
    \caption{\label{F2}
    Comparison of solutions obtained by the (mean-based) TSP and SVF algorithms, and the optimal solution. 
    The red node denotes the depot, while the blue nodes correspond to clients.
    The node attributes represent the SCV of each node, which is to be interpreted as a measure for the variability of the corresponding service time; mean service times equal~1 and travel times are assumed deterministic. The edge attributes display the optimal inter-appointment times calculated for the corresponding tour.
    }
  \label{fig:toy-example-gamma5}
\end{figure}

The intrinsic difficulty of the RAS problem lies in the fact that it is not straightforward to combine algorithms from appointment scheduling and routing (TSP), as can be seen as follows. 
A na\"{\i}ve method would be to replace the travel times and service times by their respective means, and to use off-the-shelf techniques to solve the resulting TSP (`mean based TSP', that is). 
The problem with this approach, however, is that the clients early in the tour could be the ones corresponding to a large variance, thus contradicting the principles underlying SVF.
Informally, one could say that the solution to the RAS problem finds a compromise between TSP and SVF. This effect is illustrated by Figure \ref{F2}: it shows that in this instance TSP leads to a route that differs strongly from the SVF-based one, while the optimal route (determined using full enumeration) combines elements of both.

A standard method to incorporate randomness is through the {\it sample average approximation} approach \cite{kleywegt2002,ruszczynski2003}, which involves repeatedly sampling random variables to create a large number of scenarios.
Generally, a mixed-integer linear program (MILP), which includes both routing and appointment scheduling decisions, is solved across all these scenarios.
Solving such MILPs with off-the-shelf optimization software is ineffective, but tailor-made solutions, such as the L-shaped method, make it possible to optimally solve instances with up to ten clients \cite{zhan2021}.
However, dealing with instances beyond this size remains highly challenging.

\vspace{3mm}

{\it Contributions.} 
Our approach radically differs from earlier studies. Rather than working with sampled scenarios, we incorporate the randomness into the objective function, which is subsequently optimized using soundly tuned metaheuristics. In more detail, our contributions are the following:
\begin{itemize}
    \item[$\circ$] We cast the objective function, for a given route and given appointment times, in terms of an appropriately chosen {\it queueing model}. 
    In this queueing model, the arrival times are deterministic, whereas we represent the service times by random quantities with {\it phase-type distributions}. 
    This class of distributions combines attractive properties: it is capable of closely approximating any given distribution \cite{ala2022}, but at the same time it allows a high level of mathematical tractability.
    Building on the ideas of \cite{dueck1993,kaandorp2007}, we are thus able to establish an explicit expression for our objective function.
    \item[$\circ$] Due to the large matrices involved in the phase-type method described above, the evaluation of the objective function becomes prohibitively slow for problem instances of realistic size. 
    We remedy this by applying a {\it heavy-traffic approximation}, which provides an easy-to-compute, closed-form expression for the objective function. 
    This expression allows explicit optimization over the appointment times (for a given route, that is). While the heavy-traffic approximation is not necessarily precise, it succeeds in capturing the essence of the underlying dynamics. 
    To have the best of both worlds, we use a hybrid approach for determining the optimal route in our algorithm. 
    \item[$\circ$] The optimal route is found relying on the well-developed and well-tested body of algorithms known as {\it metaheuristics} \cite{gendreau2010}. More concretely, we implement a {\it large neighborhood search} metaheuristic \cite{Pisinger2019}, which improves a tour by iteratively destroying and repairing parts of the route. 
    Combined with the fast heavy-traffic approximation, our proposed technique can handle instances of up to 40 clients within a reasonable amount of time. 
    Through comprehensive numerical experiments we show that our approach improves over other algorithms, including the modified TSP heuristic from \cite{zhan2021}.
\end{itemize}

{\it Organization.}
This paper is organized as follows. In Section \ref{lit} we provide an account of the relevant literature, both on the routing and the appointment scheduling component. Section \ref{sec:problem} formally defines the RAS problem.
Then, Section \ref{EVA} provides a closed-form expression by which the objective function can be evaluated, as well as its explicit heavy-traffic approximation. It also demonstrates that the heavy-traffic approximation captures the essence of the underlying random dynamics. The large neighborhood search metaheuristic is presented in Section \ref{sec:lns}.
Section~\ref{sec:experiments} presents a large set of numerical experiments. 
Finally, Section~\ref{sec:conclusion} concludes the paper.

\section{Literature}
\label{lit}

In this section, we consider some related literature concerning appointment scheduling (AS), routing, and the combined routing and appointment scheduling (RAS).

As mentioned, appointment scheduling is one of the two aspects of RAS problems. In particular, when we fix the route and focus on the appointment times, then our problem is reduced to AS. 
There is a long tradition of AS problems, initiated by the seminal paper of Bailey and Welch \cite{welch1952} who consider an appointment system for a hospital's outpatient department. As the performance analysis of an AS system is clearly non-trivial, some early papers have evaluated different appointment rules using simulation \cite{fetter1966, ho1992}. However, the majority of papers focus on finding optimal appointment times. 
For example, \cite{denton2003} use a sequential bounding approach to frame the problem as a linear program and solve it using an $L$-shaped method. Some other approaches to finding optimal appointment times are, for instance, simulation optimization \cite{klassen2009, klassen2019} and methods based on local search \cite{kaandorp2007, zacharias2020}. In fact, \cite{zacharias2020} show that the AS problem for discrete appointment epochs is multimodular, allowing for an efficient optimization approach. 

Another stream of research employs phase-type approximations for the service times in AS problems. Phase-type distributions can approximate any distribution arbitrarily closely and they allow for exact results using matrix algebraic manipulations; see \cite{bosch2001, wang1997} for some early studies. In the current paper, we mainly follow the approach of \cite{kuiper2015, kuiper2023} for an exact evaluation. In those studies, the authors fit a phase-type distribution to the service times based on the first two moments. The optimal appointment times are determined numerically by exploiting the convexity of the objective function \cite{kuiper2023}, which is a weighted combination of idle and waiting times in the AS setting.
Another approach to determine optimized appointment times is to use stationary models and/or heavy traffic approximations; such approximations have been applied in \cite{kuiper2017} for the single-server case and in \cite{kuiper2022} for the multi-server situation.
We refer to \cite{ahmadi2017, ala2022, cayirli2003, gupta2008} for some review papers about AS, which are all motivated by applications in health care.

The second aspect of our RAS problem is the optimization of the route. If we neglect the waiting and idle times at the nodes and only consider travel time, then the problem reduces to a TSP, which is closely related to the class of vehicle routing problems (VRPs).
Both the TSP and VRP are classical problems in combinatorial optimization; see e.g. \cite{cook2015} for a historic treatment of TSP and \cite{konstantakopoulos2020} for a recent survey on VRP.
Many heuristics and metaheuristics have been developed to solve VRPs, of which large neighborhood search (LNS) has been very succesful.
We refer to \cite{Pisinger2019} for some applications of LNS to variants of TSP and VRP.

Moreover, there is a limited amount but a variety of literature concerning combined routing and appointment scheduling. For instance, a stream of papers considers VRP models with time window assignment, in which a time window for each client needs to be determined before demand becomes known \cite{dalmeijer2018, spliet2015}. Typically, these models focus on minimizing travel and/or operator costs and impose hard time window constraints. Other studies focus on robust optimization for VRP problems with time windows under uncertainty \cite{hu2018, shi2019}. Furthermore, \cite{zhan2018} considers a stochastic programming approach in which flexible appointment times are allowed with a penalty on idle and waiting time. 
These papers incorporate uncertainty in their service times using scenarios that stem from the stochastic and robust programming domains, which differs from the way in which we incorporate uncertainty in our model. 
A difficulty with a scenario-based approach is that if the number of nodes becomes larger, it will be more challenging to represent the joint service times by representative sets of scenarios. In addition, the size of stochastic programs grows with the number of scenarios, posing computational challenges \cite{powell2019}. Instead, the basics of our approach stem from queueing-theoretic dynamics, as in AS problems.

Finally, the two papers that consider models that are most closely related to our model are \cite{tsang2023, zhan2021}. The RAS problem in \cite{zhan2021} is similar to ours, albeit that we also consider stochastic travel times. Again, the fundamental difference is the use of scenarios to represent uncertainty in service times. The authors of \cite{zhan2021} develop a MILP and an $L$-shaped method next to an easy-to-implement heuristic. In the experiments, they manage to handle instances of up to 10 nodes. The model in \cite{tsang2023} is similar, but the authors additionally consider two distributionally robust optimization models. Similar to the results presented in \cite{zhan2021}, their experiments have up to 10 clients. Our work shows that the application of appropriate queueing models within the optimization allows us to consider larger instances.


\section{Model formulation}
\label{sec:problem}

  In this section, we formulate the integrated routing and appointment scheduling (RAS) problem.
  We consider a single operator situated at a starting location denoted by the index $0$.
  The operator has to make a tour along $n$ client locations, indexed by $1, \ldots, n$, after which it returns to its starting location.
  The travel time between two locations $i, j = 0, 1, \ldots, n$ with $i \neq j$ is denoted by $T_{ij}$ and the service time at client $i$ is denoted by $B_{i}$. Both $T_{ij}$ and $B_i$ are assumed to be non-negative, independent and identically distributed random variables. A special case of specific interest is when $T_{ij}$ and $B_i$ follow a phase-type distribution.


  We define a tour $\Tour = (\tour_{1}, \ldots, \tour_{n})$ as a sequence of client indices, where $\tour_{j}$ denotes the index of the $j$-th visited client location.
  For notational convenience, we assume $\tour_0 = \tour_{n+1} = 0$, i.e., the starting location.
  For a given tour $\Tour$, the operator assigns to each client $\tour_j$ an appointment time $t_{\tour_j} \geq 0$ at which the operator plans to provide service to the client, and we assume that $t_{0} = 0$.
  The inter-appointment time $x_{\tour_{j}} \coloneqq t_{\tour_{j}}-t_{\tour_{j-1}}$ denotes the time between the appointment times of clients $\tour_{j-1}$ and $\tour_j$.
  
 When the operator arrives at a client earlier than scheduled, we say that there is \textit{idle time}, whereas if the operator arrives later than scheduled, we say there is \textit{waiting time}.
   Throughout, we denote by $I_{\tour_{j}}$ the idle time prior to the appointment time at client $\tour_j$, whereas $W_{\tour_{j}}$ is the waiting time of client $\tour_j$; see part (B) of Figure~\ref{F1} for a schematic representation.
  With slight abuse of notation, we denote by $T_{\tour_j} \coloneqq T_{\tour_{j-1}, \tour_j}$ the travel time from location $\tour_{j-1}$ to location $\tour_{j}$.
  For $j=1, \dots, n$, the following recurrence relations apply:
  \begin{align*}
    I_{\tour_j} &=\max\{x_{\tour_{j}}-T_{\tour_j}-W_{\tour_{j-1}}- B_{\tour_{j-1}},0\},\\
    W_{\tour_j} &=\max\{T_{\tour_j}+W_{\tour_{j-1}}+ B_{\tour_{j-1}}-x_{\tour_{j}},0\},
  \end{align*}
  with $W_0 = 0$ and $B_0 = 0$.
  In queueing theory, these relations are known as the {\it Lindley recursion} \cite[Section III.6]{asmussen2003}.
  They define a queueing model in which the inter-arrival times are given by the $x_{\tour_j}$ and the service requirements by the random variables $U_{\tour_j}:=T_{\tour_j}+ B_{\tour_{j-1}}$.
  Note that $U_{\tour_1} = T_{\tour_1}$, i.e., $U_{\tour_1}$ equals the travel time from the depot to the first location $\tour_1$ and does not include any service time.

Our goal is to find a tour $\Tour \equiv (\tour_{1}, \ldots, \tour_{n})$ and inter-appointment times $\IA \equiv (x_{\tour_{1}},\ldots,x_{\tour_{n}})$ that minimize an objective function that is a weighted sum of three elements: (i)~the expected value of the total travel time of the tour, (ii)~the expected aggregate idle times of the operator, and (iii)~the expected waiting times of the individual clients.

Hence, our objective function is, for weights $\omegaT, \omegaI, \omegaW_j \geq 0$,
\begin{equation}
  \label{objf}
  \objfunnew =
  \omegaT \sum_{j=1}^{n+1} \EX T_{\sigma_j}
  +\omegaI \sum_{j=1}^n \EX I_{\sigma_j} +  \sum_{j=1}^n \omegaW_{\sigma_j}\, \EX W_{\sigma_j}.
\end{equation}
Note that we allow clients to have heterogeneous waiting time weights $\omegaW_j$. 
Moreover, observe that the quantities ${\mathbb E}\,T_{\sigma_j}$ depend on the ordering $\Tour$, whereas the quantities ${\mathbb E}\,I_{\sigma_j}$ and ${\mathbb E}\,W_{\sigma_j}$ depend on both the ordering $\Tour$ and the inter-appointment times ${\bs x}$.
The RAS problem is thus defined as
\begin{equation}\label{eq:RAS}
\mathscr{L}^\star = \min_{\Tour \in \mathscr{S}, \IA \in \mathbb{R}_+^n} \objfunnew, \tag{RAS}
\end{equation}
where $\mathscr{S}$ is the collection of permutations of $\{1, \ldots, n\}$.

The RAS problem is challenging for two main reasons.
First, it requires the exploration of a complex, integrated search space of routing and appointment scheduling decisions, where the latter depends on the chosen tour.
Second, the evaluation of $\objfunnew$ for a given ordering $\Tour$ and given inter-appointment times $\IA$ requires the computation of expected values of the idle times $I_j$ and the waiting times $W_j$, which is effectively equivalent to computing their full distributions.
This is, in general, problematic: for general travel and service time distributions, one cannot obtain closed-form expressions for the density or distribution function of the idle and waiting times.
In the next section, we show that there are accurate and effective approximation procedures from queueing theory that can provide closed-form expressions.

\section{Queueing approximations for the scheduling problem}\label{EVA}

In this section, we consider the appointment scheduling part of the RAS problem. 
We present two queueing-based approximation procedures to evaluate the objective function $\objfunnew$ for a given tour $\Tour$ and inter-appointment times ${\bs x}$, as well as the objective function that is optimized over the inter-appointment times $\min_{{\bs x\in {\mathbb R}_+^n}} {\mathscr L}(\Tour,\IA)$. In Subsection~\ref{sec:algo_objective} we give a closed-form expression for the objective function exploiting the structure of phase-type distributions.  
Then, in Subsection~\ref{sec:heavy-traffic}, we propose an approximation of the objective function based on a heavy-traffic regime, leading to an expression that allows an explicit minimization over ${\bs x\in {\mathbb R}_+^n}$. A hybrid version that can efficiently compute good appointment schedules is introduced in Subsection~\ref{sec:hybrid}.

Throughout this section, we assume that the tour $\Tour$ is given. To avoid excessive notation, we let $j$ here denote the index of the $j$-th visited client in the given $\Tour$ (that is, we use $j$ instead of $\tour_j$). 
Furthermore, the analysis here focuses on the idle and waiting times, i.e., the second and third terms of the objective function in (\ref{objf}).

\subsection{Phase-type method to evaluate objective function}\label{sec:algo_objective}

As alluded to in Section~\ref{sec:problem}, the evaluation of $\objfunnew$ is challenging, even for a given tour $\Tour$ and given inter-appointment times $\IA$. This can be observed from the Lindley recursion; the evaluation of ${\mathbb E}\,W_j$ requires that we know the distribution of $W_{j-1}$.
One cannot obtain closed-form expressions for the distribution function of the idle and waiting times for general travel and/or service times.
This motivates why we follow a phase-type method for approximating ${\mathbb E}\,I_j$ and ${\mathbb E}\,W_j$ for any travel and service time distribution, which has been advocated, and firmly backed, in \cite{kuiper2023}. 
Before presenting the two-step method, let us first give some relevant preliminaries and motivation related to phase-type distributions. 

We first recall that any distribution on the positive half line can be approximated arbitrarily closely by a {\it phase-type distribution}; we formally say that the class of phase-type distributions is dense in the class of all distributions on $(0,\infty)$ \cite[Thm.\ III.4.2]{asmussen2003}. The class of phase-type distributions contains mixtures and
convolutions of exponential distributions, and thus includes (mixtures of) Erlang distributions and hyperexponential distributions; formally a phase-type random variable is defined as the entrance time of an absorbing state in a continuous-time Markov chain \cite[Section III.4]{asmussen2003}. 
In particular, a phase-type distribution is characterized through the triple $(d, {\bs\alpha}, V)$ with $d$ the number of transient states of the Markov chain, ${\bs\alpha}$ an initial distribution, and $V$ the transition rate matrix. 
We refer to Appendix \ref{APP1} for some background on phase-type distributions; in particular, we provide the specific definitions of the ones that play a role in this paper. 

In \cite{wang1997} it is pointed out how the distributions of $I_j$ and $W_j$ can be (numerically) evaluated in case $U_j=T_j+B_{j-1}$ is of phase-type; indeed, applying techniques from linear algebra, a recursive technique to determine the cumulative distribution functions ${\mathbb P}(I_j\leqslant t)$ and ${\mathbb P}(W_j\leqslant t)$, for $t\geqslant 0$, has been devised. We adopt a slightly different approach to evaluate the objective function $\objfunnew$.

Now, our phase-type method consists of the following two steps: (1) approximate the sum of travel and service times for each client by a phase-type distribution, and (2) provide an explicit formula for the objective function in case of phase-type service and travel times. 
For step~(1), the specific phase-type distributions that we will be using are mixtures of Erlang distributions and the hyperexponential distribution, see \cite{tijms1986} for additional justification. 
The underlying idea is to fit the parameters of these distributions such that they match the first two moments of the distribution of the target random variable, say $U$ (or, equivalently, the first moment ${\mathbb E}\,U$ and the {\it squared coefficient of variation} $\mbox{\scv}(U):= {\mathbb V}{\rm ar} \,U/({\mathbb E}\, U)^2$). Fitting the first two moments means that evidently some detail about the distribution of $U$ is lost. It should, however, be borne in mind that (i) in general only lower moments of the random variables under consideration can be reliably estimated from historical data, and (ii), as backed by the experiments in \cite{kuiper2016}, the approximation of random variables by their phase-type counterparts is highly accurate. In the sequel we denote by $\bar U$ the phase-type random variable that approximates $U$.

\begin{itemize}
  \item In the case of $\bar U$ being a mixture of Erlang distributions,  $\bar U$ corresponds to an Erlang distributed random variable with $K - 1$ phases and mean $(K - 1)/\mu$ with probability $p$, and with an Erlang distributed random variable with $K$ phases and mean $K/\mu$ with probability $1 - p$; we write that $\bar U$ is distributed as $E_K(\mu,p)$. As pointed out in \cite{tijms1986}, for such a $\bar U$, with $K\in\{2,3,\ldots\}$, the squared coefficient of variation (SCV) satisfies:
\[\mbox{\scv}(\bar U) = \frac{{\mathbb V}{\rm ar} \,\bar U}{({\mathbb E}\,\bar U)^2}= \frac{K-p^2}{(K-p)^2}\in \left( K^{-1},{(K-1)}^{-1} \right] .\]
        As a consequence, if $\mbox{\scv}(U)\in \left(K^{-1},(K-1)^{-1}\right]$, then we propose to approximate $U$ by $E_K(\mu,p)$ with
    \[ p = \frac{1}{1+\mbox{\scv}(U)} \left( \mbox{\scv}(U) K - \sqrt{(1+\mbox{\scv}(U)) K - \mbox{\scv}(U) K^2} \right)
    \quad {\rm and} \quad
    \mu = \frac{K-p}{{\mathbb E}\,U},
    \]
        such that $U$ and $\bar U$ have the same first two moments, cf.~\cite{tijms1986}.
 \item In the case of $\bar U$ being hyperexponential, $\bar U$ corresponds with probability $p$ to an exponential distribution with rate $\mu_1$, and with probability $1-p$ to an exponential distribution with rate $\mu_2$; we write that $\bar U$ is distributed as $H_2(\mu_1,\mu_2,p)$.
 When fitting the first two moments with three parameters, there is clearly a degree of freedom. To remedy this, we choose $\mu_1 =2p\mu$ and $\mu_2 =2(1-p)\mu$ for some $\mu >0$ (a choice referred to as `balanced means'). An elementary computation reveals that
 \[\mbox{\scv}(\bar U) = \frac{1}{2p(1-p)},\]
  which is larger than or equal to 1. 
  Hence, if $\mbox{\scv}(U) \ge 1$ it can be verified (see \cite{tijms1986}) that one can pick
  \[p= \frac{1}{2}
  \left(1 + \sqrt{{\frac{\mbox{\scv}(U)-1}{\mbox{\scv}(U)+1}}}\right),
 \quad \mu_1 = \frac{2p}{{\mathbb E}\,U}, 
 \quad \mu_2 = \frac{2(1-p)}{{\mathbb E}\,U},
  \]
    to match the first two moments of $U$ and $\bar U$. 
\end{itemize}

Turning to step~(2), we present an approach in Appendix~\ref{APP2}, based on \cite{wang1997}, to evaluate the objective function $\objfunnew$ using that all $\bar U_j$ are of phase type. More specifically, we identify the distribution of the random variables $R_j:=W_j+\bar U_{j+1}$, which is in fact of phase type itself. As shown in the appendix, the phase-type distribution of $R_j$ is characterized by the triplet $(D_j, {\bs\alpha}_{(j)}, V^{(j)})$; the parameters are recursively defined and can be found in Appendix~\ref{APP2}.
An explicit formula for the objective function is given in the following theorem, whereas the derivation is deferred to Appendix~\ref{APP2}. Observe that the expression for the objective function is exact in case all $U_j$ are phase type; if not all $U_j$ are phase-type, the theorem provides an accurate approximation.

\begin{theorem}\label{thm:obj}  If all $U_j$ are phase type, then, for $\Tour \in{\mathscr I},{\bs x\in {\mathbb R}_+^n}$,
\begin{align*}{\mathscr L}&(\Tour,{\bs x)} = \omegaT
\sum_{j=1}^{n+1} {\mathbb E}\,T_{j} \:+\\& \sum_{j=1}^n \Bigg( \omegaI \left(x_{j}+{\bs\alpha}_{(j)}(V^{(j)})^{-1} {\bs 1} - {\bs\alpha}_{(j)}(V^{(j)})^{-1} e^{V^{(j)}x_{j}}{\bs 1}\right) 
+ \omegaW_j \left(-{\bs\alpha}_{(j)}(V^{(j)})^{-1} e^{V^{(j)}x_{j}}{\bs 1}\right) \Bigg),\end{align*}
with ${\bs\alpha}_{(j)}$ and $V^{(j)}$ recursively defined in Appendix~\ref{APP2} and ${\bs 1}$ a vector of ones.
\end{theorem}

When it comes to our RAS problem, the above methodology has one serious limitation.
The expressions in Appendix~\ref{APP2} reveal that, in particular when some of the random variables $U_j$ have a low SCV, the dimension of the matrices can become prohibitively large.
In addition, evaluating the objective function requires a number of matrix inverses and matrix exponentials, which are computationally expensive.
The optimal appointment times can be numerically obtained, e.g., using gradient descent methods, by exploiting the fact that the objective function is convex in ${\bs x}$, see \cite{kuiper2023}.
Although a single evaluation does not take much time in practice, obtaining the optimal solution 
\[{\mathscr L}_{\rm ph}(\Tour) := \min_{{\bs x\in {\mathbb R}_+^n}} {\mathscr L}(\Tour,{\bs x)}, \]
requires many iterations, for which this approach becomes infeasible.
Instead of optimizing over the inter-appointment times, we present a heavy-traffic approximation in the next section.

\subsection{Heavy-traffic approximation}
\label{sec:heavy-traffic}

In this subsection we discuss a heavy-traffic approximation of our objective function, which allows us to easily optimize over the inter-appointment times.
It is expected to be particularly accurate when one cares more about idle times than about waiting times.
The reason is that in that regime the schedule will be such that the queueing system, as defined through the Lindley recursion, will be operating under heavy load.

The first step is that, based on the heavy-traffic approximation for the mean waiting time in a GI/G/1 queue
\cite[Section X.7]{asmussen2003}, we propose the approximation
\[{\mathbb E}\, W_j \approx \frac{S_{j}}{2(x_j-\mathbb {E} \,U_j)},\]
where we denote 
\[S_j:=\frac{\sum_{i=1}^j \beta^{j-i} {\mathbb V}{\rm ar}\,U_i}{\sum_{i=1}^j \beta^{j-i}},  \]
where empty sums being defined to equal $0$.
The constant $\beta \in [0,1]$ represents the rate of the exponential decay in the variance of the service and travel time that a client has on subsequent clients. 
After some numerical experiments, we chose $\beta = 0.5$.
As this approximation is based on results for a stationary random variable, we assume that $x_j > \mathbb {E} \,U_j$.
In addition, for the mean idle times, we can use the intuitive approximation
\[{\mathbb E}\,I_j \approx x_j - \mathbb {E} \,U_j.\]
Upon combining the above approximations, we obtain the following heavy-traffic approximation of the objective function:
\[ {\mathscr L}_{{\rm ht}}(\Tour,{\bs x}):=\omegaT
\sum_{j=1}^{n+1} {\mathbb E}\,T_{j}+\omegaI\sum_{j=1}^n( x_j - \mathbb {E} \,U_j)  +  \sum_{j=1}^n \omegaW_{j} \frac{S_{j}}{2(x_j - \mathbb {E} \,U_j)}.\]
We can thus try to find an approximate value of ${\mathscr L}^\star$ by performing the optimization problem
\[{\mathscr L}^\star_{\rm ht} := \min_{\Tour \in{\mathscr I},{\bs x\in {\mathbb R}_+^n}} {\mathscr L}_{\rm ht}(\Tour,{\bs x)}.\tag{RAS-ht}\]
Interestingly, it is possible to {\it explicitly} solve the optimization over ${\bs x}$.
To this end, first note that it is straightforward to verify that this function is convex in its arguments $x_1,\ldots, x_n$. By solving its first-order conditions, we have that
for $j=1,\ldots, n$,
\begin{equation}
    \label{htopt1}
x_j^{\rm (ht)} = {\mathbb E}\,U_{j} + \sqrt{\frac{\omegaW_j S_j}{2 \omegaI}}. \end{equation}
These inter-appointment times are intuitively appealing and turn out to provide good appointment schedules.
Note that these findings are consistent with those for the stationary homogeneous case under heavy traffic, as have been presented in \cite{kuiper2017}. 
As an aside, we mention that the corresponding value of the objective function can now be evaluated in closed form.

\begin{proposition}\label{PHT}
For all $\Tour \in{\mathscr I}$,
\begin{equation}
    \label{COST1}
{\mathscr L}_{{\rm ht}}(\Tour):=\min_{{\bs x\in {\mathbb R}_+^n}} {\mathscr L}_{\rm ht}(\Tour,{\bs x)}={\mathscr L}_{\rm ht}(\Tour,{\bs x}^{\rm (ht)})=\omegaT
\sum_{j=1}^{n+1} {\mathbb E}\,T_{j}+\sqrt{2\,\omegaI}\, \sum_{j=1}^n\sqrt{\omegaW_j \, S_{j}},\end{equation}
with ${\bs x}^{({\rm ht})}$ given by \eqref{htopt1}.
\end{proposition}

\begin{remark} {\em
Importantly, for the case without travel times (e.g., the conventional appointment scheduling situation where ${\mathbb E}\,T_j=0$) the expression \eqref{COST1} for the approximate minimum cost offers additional support for the SVF rule: as  can be seen directly, it is minimized if the clients are sequenced in increasing order of variance. 
}
\end{remark}

\subsection{Hybrid approximation} 
\label{sec:hybrid}


In Subsection~\ref{sec:algo_objective}, we presented a phase-type method to obtain an accurate approximation of the objective function ${\mathscr L}(\Tour,{\bs x})$; in fact, this method is exact when travel and service times follow a phase-type distribution. Although this method is accurate, it is relatively slow, in particular if we aim to find optimized inter-appointment times ${\bs x}$. In Subsection~\ref{sec:heavy-traffic}, we derived a heavy-traffic approximation, which is less accurate than the phase-type method, but it is fast. 
Based on numerical experiments, presented in Appendix~\ref{accht}, we propose a hybrid version that combines the strengths of the two approaches:
\begin{equation}\label{defhyb}
{\mathscr L}_{{\rm hyb}}(\Tour):= {\mathscr L}(\Tour,{\bs x}^{\rm (ht)}),
\end{equation}
with ${\bs x}^{({\rm ht})}$ the vector whose entries are given by \eqref{htopt1}.
This hybrid version can be seen as `best of both worlds': the appointment times stem from the heavy-traffic approximation, i.e., require virtually no computation time. These appointment times are then inserted into the phase-type based objective function presented in Theorem~\ref{thm:obj} and yields a more accurate proxy of the objective function than when ${\bs x}^{\rm (ht)}$ is inserted into the heavy-traffic based objective function, i.e., Proposition~\ref{PHT}.



\section{Large neighborhood search metaheuristic}
\label{sec:lns}
The previous section presented an approximation to the RAS problem, transforming it to a pure routing problem.
In this section, we propose the use of a large neighborhood search (LNS) metaheuristic to efficiently explore the search space of routing solutions.
LNS was first introduced by \cite{Shaw1998}; a recent and more general treatment is given by \cite{Pisinger2019}, including many successful examples of the LNS metaheuristic in solving combinatorial optimization problems.
The main idea underlying LNS is to improve an initial solution (i.e., a tour) by iteratively \textit{destroying} and \textit{repairing} the solution.

\subsection{Algorithm outline}
A general outline of the LNS metaheuristic is given in Algorithm~\ref{alg:lns}. 
The algorithm takes an initial solution and initializes it as the best solution.
In each iteration of the algorithm, a random pair of destroy $d(\cdot)$ and repair $r(\cdot)$ operators is selected and applied to the current tour, resulting in a candidate solution $\Tour^c$.
An acceptance criterion determines whether this candidate solution is accepted as the solution for the next iteration or rejected in favor of the old solution $\Tour$.
Moreover, if the identified candidate solution improves over the best-found solution so far, we register it as the new best solution $\Tour^b$.
We iterate this procedure until a provided time limit is met, after which we return the best-found solution. 

\begin{algorithm2e}
\caption{Large neighborhood search}
\label{alg:lns}
\DontPrintSemicolon
\KwIn{Initial solution $\Tour$}
\KwOut{The best-found solution $\Tour^b$}
$\Tour^b \gets \Tour$\;
\While{\text{time limit is not exceeded}}{
    Select destroy and repair operator pair $(d, r)$ at random\;
    $\Tour^c \gets r(d(\Tour))$\;
    \If{$\Tour^c$ is accepted}{
        $\Tour \gets \Tour^c$\;
    }
    \If{$\Tour^c$ has a better objective value than $\Tour^b$}{
        $\Tour^b \gets \Tour^c$\;
    }
}
\Return{$\Tour^b$}\;
\end{algorithm2e}

\subsection{Details of the algorithm}
This section describes the details of the LNS algorithm in the context of the RAS problem. 

\subsubsection{Objective function}
A candidate solution $\Tour^{c}$ is evaluated using the objective function $\objfunhyb(\Tour^{c})$ defined in (\ref{defhyb}).

\subsubsection{Initial solution}
We use a random permutation of client visits as the initial solution.

\subsubsection{Destroy and repair operators}
We consider two destroy operators. 
The first one is the \textit{random client} operator, which randomly removes between 1 and $D$ clients from the current solution, with $D\in{\mathbb N}$ a parameter of the algorithm.
The second one is an \textit{adjacent client} operator, which removes a randomly selected sequence of between 1 and $D$ adjacent clients.
For the repair operator, we consider a \textit{greedy insert} repair operator.
This operator inserts each of the removed clients into the position that increases the objective function by the least amount, one by one.

\subsubsection{Acceptance criterion}
We use the record-to-record-travel acceptance criterion \cite{dueck1993}.
This criterion accepts a new candidate solution $\Tour^c$ if $\mathscr{L}_{\rm hyb}(\Tour^c) - \mathscr{L}_{\rm hyb}(\Tour^b)$ is smaller than some updating threshold $H$ for the current best-found solution $\Tour^b$.
This threshold is initialized at some starting value $H_0 \geq 0$ and is equal to $H = H_0 \cdot t/t_{\max}$, where $t$ denotes the elapsed runtime and $t_{\max} \geq t \geq 0$ denotes the maximum time limit.

\subsubsection{Post-termination optimization}
\label{sec:post-termination}
After termination of the LNS algorithm, we use the best-found tour $\Tour^b$ to compute the optimal appointment schedule $\IA^\star$, i.e.,
\begin{equation}
    \IA^\star = \argmin_{\IA \in \mathbb{R}^n_+} \mathscr{L}(\Tour^b, \IA).
\end{equation}
The final solution to the RAS problem is then given by $(\Tour^b, \IA^\star)$.
In summary, the idea is that to determine the tour, we solely use the hybrid approximations with appointment schedules based on heavy traffic; for the best-found tour, we finally determine the appointment schedule based on the phase-type method of the objective function.

\section{Numerical experiments}
\label{sec:experiments}
In this section, we present the numerical experiments. 
We describe the generation procedure of benchmark instances in Subsection~\ref{sec:benchmark-instances}, the implemented benchmarks algorithms in Subsection~\ref{sec:benchmark-algorithms}, parameter tuning in Subsection~\ref{sec:parameter-tuning} and the results in Subsection~\ref{sec:results}.
All instances, algorithms, and results are openly available at \url{https://github.com/leonlan/routing-appointment-scheduling}.

\subsection{Instance generation procedure}\label{sec:benchmark-instances}
We followed the instance generation procedure described in \cite{zhan2021}, which is common in the integrated routing and appointment scheduling literature.
An instance is defined on a $\left[0, 50\right]^2$ square grid in which the depot is located at coordinate $(0, 0)$ and $n$ clients are uniformly randomly dispersed.
The mean travel time between two locations equals the Euclidean distance, and we consider a fixed SCV for travel times of $0.15$.
The mean service times are sampled from $U(30, 60)$, where $U(a,b)$ denotes the uniform distribution on $[a,b]$.
The SCVs for service times are either all sampled from $U(0.15, 0.5)$ or from $U(0.5, 1.5)$. We refer to these two cases as \textit{low} and \textit{high} service time SCVs.
The weights for client waiting times $\omegaW_{j}$ are sampled from $U(1, 10)$, and the weight of the operator’s idle time $\omegaI$ is fixed at 2.5.
To investigate the trade-off between travel times and appointment costs, we consider three scenarios with $\omegaT$ being 0.5, 1, and 2, respectively.

For a given instance size $n$, we sampled 20 different realizations for all parameters except for the SCVs of service times and the weights for travel times.
Using these base instances, we sampled the SCVs of service times for the low and high variance scenarios and used each of the three different $\omegaT$ values.
This resulted in $20 \times 2 \times 3 = 120$ unique problem instances for a given instance size $n$.

\subsection{Benchmark algorithms}\label{sec:benchmark-algorithms}
We implemented three benchmark heuristics to compare against our LNS metaheuristic.
All three heuristics first identify a tour of the client visits, after which the optimal appointment schedule is computed, similar to the procedure outlined in Subsection~\ref{sec:post-termination}.
\begin{itemize}
    \item[$\circ$] The first heuristic is the modified TSP (MTSP) heuristic proposed by \cite{zhan2021}. The MTSP solves a deterministic TSP in which the edge weights are modified to take into account appointment scheduling costs.
A more detailed explanation of the algorithm is given in Appendix~\ref{APP2}.
\item[$\circ$] The second heuristic solves the deterministic TSP using the mean travel times, yielding the tour
    $\Tour^{\rm TSP} = (\tour_1^{\rm TSP}, \ldots, \tour_n^{\rm TSP})$.
        For this tour, the objective value follows by optimizing over the appointment times for tour $\Tour^{\rm TSP}$, i.e., $\min_{{\bs x\in {\mathbb R}_+^n}} {\mathscr L}(\Tour^{\rm TSP},\IA)$.
        Since the benchmark instances assume Euclidean distances, the reverse direction of the tour $(\tour_n^{\rm TSP}, \ldots, \tour_1^{\rm TSP})$ provides the same travel time.
        Hence, we take the direction of the tour with the smallest objective value of the two directions.
        From our numerical experiments, we have observed that a considerable gain can be achieved by considering the two-sided tour compared to the one-sided version.
\item[$\circ$] The third heuristic is a modified version of the smallest variance first rule (MSVF) that incorporates stochastic travel times.
Starting at the depot, we select the next, not-yet-visited client with the lowest variance of combined travel and service time.
That is, given a partial tour $\Tour'= (\tour_{1}, \dots, \tour_{k})$, we select the next location by $$\argmin_{j \in V, j \notin \Tour'} {\mathbb V}{\rm ar}(T_{\tour_{k}, j} + B_j).$$
We repeat this procedure until all clients have been visited.
For the identified tour, we compute the optimal inter-appointment times similar to the TSP algorithm as described above.
When the variance of the travel times goes to zero, MSVF reduces to the classical SVF rule where clients are visited in order of smallest variance of the service times.
\end{itemize}

\subsection{Parameter tuning}\label{sec:parameter-tuning}
For LNS, we performed a full factorial design using the parameters $D$ (maximum number of removed clients) and $H_0$ (starting threshold of acceptance criterion).
We restricted the values to $D \in \{2, 3, \dots 12\}$ and $H^{\text{init}} \in \{0.01, 0.025, 0.05, 0.075, 0.10\}$, where we set $H_0$ equal to $H^{\text{init}}$ multiplied with the objective value of the initial solution. 
This resulted in 55 different parameter combinations. 
All combinations were evaluated on a set of 60 training instances of size $n \in \{10,15,20\}$ generated following the same procedure as the benchmark instances.
The time limits were set to $15, 30$ and $60$ seconds for $n=10,15$ and $20$, respectively.
The best results were obtained for $D=6$ and $H^{\text{init}}=0.05$, hence we selected these parameters.

\subsection{Results}\label{sec:results}
We now present the results of the numerical experiments. 
Subsection~\ref{sec:results-small-instances} presents the results for small benchmark instances with $n \in \{6, 8, 10\}$ and Subsection~\ref{sec:results-large-instances} presents those for large instances with $n \in \{15, 20, 25, 30, 35, 40\}$.
Subsection~\ref{sec:results-runtimes} investigates the runtime performance of the LNS algorithm and of computing the optimal appointment schedules.
Finally, Subsection~\ref{sec:results-scaling-travel-times} performs a sensitivity analysis of the algorithmic performance when shifting the emphasis from the routing to the appointment scheduling component by scaling the travel times.

We implemented all our numerical experiments in Python 3.11.
We used the LKH-3 solver \cite{helsgaun2017} to obtain solutions to the deterministic TSP, and we implemented the LNS metaheuristic using the ALNS Python package \cite{wouda2023alns}.
The optimal appointment schedules were obtained using trust-region methods from the SciPy Python package \cite{virtanen2020}.
Each instance was solved on a single core of an AMD EPYC 7H12 CPU.

For each solved instance, we computed the \textit{gap}, which is defined as the relative difference in objective value between the obtained solution and the best-found solution.
The tables report the average gap over all instances for specific instance categories.

\subsubsection{Small instances}\label{sec:results-small-instances}
\renewcommand*{\arraystretch}{1.2}
\begin{table}[!ht]
\centering
\begin{tabular}{cclcccclcccc}
\toprule
 &  &  & \multicolumn{4}{c}{Low SCV} & & \multicolumn{4}{c}{High SCV} \\
 \cline{4-7} \cline{9-12}
\rule{0pt}{3ex} $n$ & $\omegaT$ &  & LNS & MTSP & TSP & MSVF &  & LNS & MTSP & TSP & MSVF \\
\midrule
\multirow[t]{3}{*}{6} & 0.5 &  & \textbf{0.80} & 5.70 & 4.50 & 11.16 &  & \textbf{2.59} & 10.96 & 6.83 & 11.00 \\
 & 1.0 &  & \textbf{0.73} & 4.82 & 3.43 & 13.16 &  & \textbf{1.96 }& 8.15 & 5.65 & 12.70 \\
 & 2.0 &  & \textbf{0.79} & 2.71 & 2.45 & 16.50 &  & \textbf{1.97} & 7.23 & 4.20 & 15.75 \\
\midrule
\multirow[t]{3}{*}{8} & 0.5 &  & \textbf{1.12} & 5.25 & 3.73 & 9.46 &  & \textbf{3.93} & 10.35 & 7.28 & 9.69 \\
 & 1.0 &  & \textbf{0.84} & 4.47 & 2.88 & 12.17 &  & \textbf{2.84} & 7.90 & 5.68 & 11.39 \\
 & 2.0 &  & \textbf{0.56 }& 3.99 & 1.96 & 16.43 &  & \textbf{2.24} & 6.36 & 3.96 & 14.96 \\
 \midrule
\multirow[t]{3}{*}{10} & 0.5 &  & \textbf{2.04} & 7.72 & 5.50 & 10.27 &  & \textbf{3.32 }& 10.47 & 8.64 & 11.16 \\
 & 1.0 &  & \textbf{1.14} & 5.86 & 4.10 & 13.71 &  & \textbf{3.00} & 8.55 & 6.82 & 13.64 \\
 & 2.0 &  & \textbf{0.88} & 4.63 & 2.80 & 19.74 &  & \textbf{1.70} & 5.82 & 4.62 & 18.56 \\
 \midrule
Avg. &  && \textbf{0.99} & 5.02 & 3.48 & 13.62 && \textbf{2.62} & 8.42 & 5.96 & 13.21 \\
\bottomrule
\end{tabular}
\caption{
Results on small benchmark instances.
Values represent the average gaps per instance size, travel time weight and SCV variation.
The best average gap per category is marked in bold.
}
\label{tab:benchmark-small}
\end{table}
The first set of numerical experiments compares the performance of LNS and the benchmark algorithms against the (near-)optimal solution on small problem instances with $n \in \{6, 8, 10\}$.
We can obtain the optimal solution by fully enumerating all possible tours and computing the corresponding optimal appointment schedules.
For instances with $n=6$ and $n=8$, the optimal solution can be computed within reasonable time limits (below, say, 1 hour).
However, for instances with $n=10$ clients, the time needed to evaluate a single instance would take roughly 500 hours on a single core; hence, for this value of $n$ we resorted to a variant of the full enumeration approach.
More, specifically, we first enumerated all $10!$ tours and evaluated them using the heavy traffic appointment schedule ${\bs x}^{({\rm ht})}$, giving ${\mathscr L}_{\rm ht}(\Tour)$ for all $\Tour \in \mathscr{S}$.
We then considered the $8!$ solutions with the lowest objective value and evaluated the tours from these solutions again using the optimal appointment schedule, yielding ${\mathscr L}_{\rm ph}(\Tour)$ for the $8!$ considered tours $\Tour$.
Although this does not guarantee an optimal solution, we believe that the obtained solutions are near-optimal and serve as reliable reference solutions.

The time limit for LNS was set to 5, 10, and 15 seconds for $n=6, 8$, and $10$, respectively.
The runtimes of MTSP, TSP, and MSVF were negligible.
Table~\ref{tab:benchmark-small} shows the results for different instance sizes $n$, travel time weights $\omegaT$, and SCV variation.
On average, LNS achieved the lowest average gap of 1.80\% on all small instances (0.99\% for low SCV and 2.62\% for high SCV), whereas the average gaps of MTSP, TSP, and MSVF were 6.72\%, 4.72\%, and 13.41\%, respectively.
The results reveal that the integration of routing and appointment scheduling decisions in our LNS metaheuristic has clear benefits over simpler heuristics.
TSP performed surprisingly well: in our results, the TSP algorithm outperformed the MTSP and MSVF heuristics.
We noticed in our experiments that evaluating both orientations for the TSP algorithm was a prominent reason for this: evaluating only one orientation of the TSP tour resulted in considerably worse solutions.
When evaluating only one orientation of the TSP tour, the results of MTSP, TSP, and MSVF were consistent with those presented in \cite{zhan2021}.
The performance of MSVF is rather poor across all instances, as the emphasis on travel times in the instances is relatively large.

When comparing results between the different values of the travel time weight, we observe that LNS obtains lower gaps for high values, that is, when there is more emphasis on routing.
Similarly, the results show that for low service time SCVs, the average gap achieved by LNS is lower than in the case of high service time SCVs (0.99\% versus 2.62\%). 
This clearly indicates that problems that have more emphasis on the appointment scheduling problem, i.e., low travel time weights and higher variance in service times, result in generally more challenging instances, clearly showing that there is still room for improvement.

\subsubsection{Large instances} \label{sec:results-large-instances}
In the second set of experiments, we compare the performance of LNS and the benchmark algorithms on larger problem sizes with $n \in \{15, 20, 25, 30, 35, 40\}$.
The time limit for LNS was set to $30, 60, 120, 240, 480$, and $960$ seconds, respectively.
Note that it is no longer feasible to compute the (near-)optimal solution using full enumeration, hence the presented gaps are with respect to the best-found solution.

\renewcommand*{\arraystretch}{1.2}
\begin{table}[!ht]
\centering
\begin{tabular}{cclcccclcccc}
\toprule
 &  &  & \multicolumn{4}{c}{Low SCV} & & \multicolumn{4}{c}{High SCV} \\
 \cline{4-7} \cline{9-12}
\rule{0pt}{3ex} $n$ & $\omegaT$ &  & LNS & MTSP & TSP & MSVF &  & LNS & MTSP & TSP & MSVF \\
\midrule
\multirow[t]{3}{*}{15} & 0.5 &  & \textbf{0.16} & 3.75 & 2.85 & 8.95 &  & \textbf{0.16} & 5.76 & 4.10 & 6.41 \\
 & 1.0 &  & \textbf{0.16} & 2.45 & 2.03 & 13.59 &  & \textbf{0.09} & 4.04 & 3.26 & 10.13 \\
 & 2.0 &  & \textbf{0.70 }& 1.88 & 1.70 & 22.15 &  & \textbf{0.18} & 2.79 & 2.19 & 17.01 \\
 \midrule
\multirow[t]{3}{*}{20} & 0.5 &  & \textbf{0.28} & 3.41 & 2.29 & 8.53 &  &\textbf{0.03} & 5.99 & 5.20 & 7.21 \\
 & 1.0 &  & \textbf{0.53} & 2.85 & 1.88 & 13.83 &  & \textbf{0.20} & 3.74 & 3.76 & 10.57 \\
 & 2.0 &  & 1.44 & 1.69 & \textbf{1.36} & 22.91 &  & \textbf{0.21} & 2.76 & 2.78 & 18.15 \\
 \midrule
\multirow[t]{3}{*}{25} & 0.5 &  & \textbf{0.30} & 2.13 & 1.67 & 7.98 &  & \textbf{0.00} & 5.96 & 5.41 & 7.26 \\
 & 1.0 &  & 1.01 &\textbf{1.00 }& 1.07 & 13.12 &  & \textbf{0.03} & 4.80 & 4.61 & 11.34 \\
 & 2.0 &  & 1.34 & 1.22 & \textbf{0.58} & 22.48 &  & \textbf{0.10} & 3.33 & 2.98 & 18.43 \\
 \midrule
\multirow[t]{3}{*}{30} & 0.5 &  & \textbf{0.21} & 3.11 & 1.64 & 7.11 &  & \textbf{0.00 }& 5.82 & 5.40 & 5.98 \\
 & 1.0 &  & \textbf{0.70} & 2.09 & 0.84 & 12.22 &  &\textbf{0.00} & 5.55 & 4.42 & 9.68 \\
 & 2.0 &  & 3.49 & 1.05 & \textbf{0.27} & 22.01 &  & \textbf{0.01} & 3.67 & 2.47 & 16.15 \\
 \midrule
\multirow[t]{3}{*}{35} & 0.5 &  & \textbf{0.61} & 1.81 & 0.86 & 6.64 &  & \textbf{0.01} & 6.19 & 5.26 & 6.57 \\
 & 1.0 &  & 1.95 & 1.39 & \textbf{0.47} & 11.92 &  & \textbf{0.11} & 4.87 & 4.14 & 10.21 \\
 & 2.0 &  & 4.97 & 1.25 & \textbf{0.17} & 21.68 &  & \textbf{0.34 }& 3.65 & 2.85 & 17.66 \\
 \midrule
\multirow[t]{3}{*}{40} & 0.5 &  & 0.57 & 1.80 & \textbf{0.46} & 5.99 &  & \textbf{0.00} & 7.19 & 6.10 & 6.68 \\
 & 1.0 &  & 1.77 & 1.46 & \textbf{0.11} & 11.31 &  & \textbf{0.00} & 6.18 & 4.98 & 10.26 \\
 & 2.0 &  & 5.85 & 1.04 & \textbf{0.06} & 21.48 &  & \textbf{0.00} & 4.15 & 2.81 & 16.65 \\
 \midrule
Avg. &  && 1.45 & 1.97 & \textbf{1.13} & 14.11 && \textbf{0.08} & 4.80 & 4.04 & 11.46 \\
\bottomrule
\end{tabular}
\caption{
Results on large benchmark instances.
Values represent the average gaps per instance size, travel time weight and SCV variation.
The best average gap per category is marked in bold.
}
\label{tab:benchmark-large}
\end{table}




Table~\ref{tab:benchmark-large} shows results for larger instances.
LNS achieved an average gap on all instances of 0.77\%, whereas MTSP, TSP and MSVF achieved an average gap of 3.38\%, 2.58\%, 12.78\%, respectively.
Although the average gap between LNS and other algorithms has decreased compared to the results on small instances, the results still show that LNS outperformed the other heuristics on average.
Specifically, for instances with high SCVs, LNS consistently outperforms the other heuristics, with an average gap difference of 3.96\% when compared to TSP.
However, LNS no longer obtains the best results for instances with low SCVs: it is generally outperformed by TSP when instances become larger and the travel time weights are higher.
This is a result of the high computational complexity of evaluating the objective function: phase-type distributions with low SCVs require more phases, and thus yield larger matrices.
We note that the performance of LNS can be improved by increasing the time limit. The impact of the runtime is discussed in more detail in the next section.

\subsubsection{Analysis of the runtimes}
\label{sec:results-runtimes}
We provide here more details about the runtimes of the algorithms.
First, we show how much time one single LNS iteration needs on average.
And, second, we show how much time is required to compute the optimal appointment schedule when the final tour has been selected.
Both evaluations rely on evaluating the objective function presented in Theorem~\ref{thm:obj}. 
The time complexity to evaluate this function grows fast: it requires general matrix operations such as multiplication and inversion, as well as the computation of the matrix exponential.
One important factor that determines the size of the matrices is the SCVs of the combined service and travel times.

\begin{figure*}[!ht]
    \centering
    \subfloat[LNS iteration]{\includegraphics[width=0.48\linewidth]{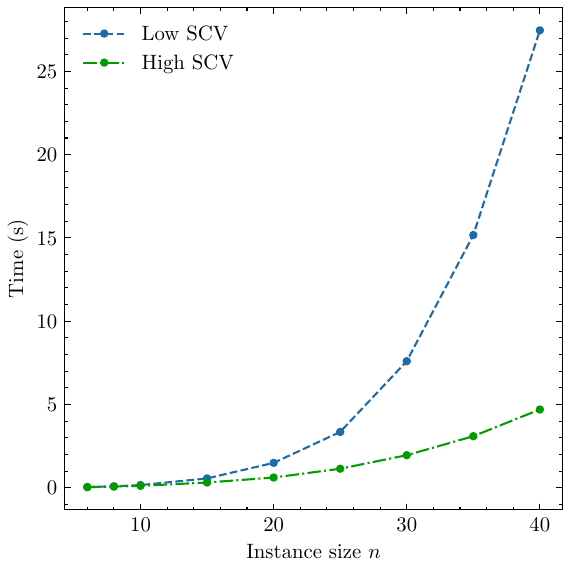}%
    \label{fig:results-runtime-lns-iteration}}
    \hfill
    \subfloat[Optimal schedule]{\includegraphics[width=0.48\linewidth]{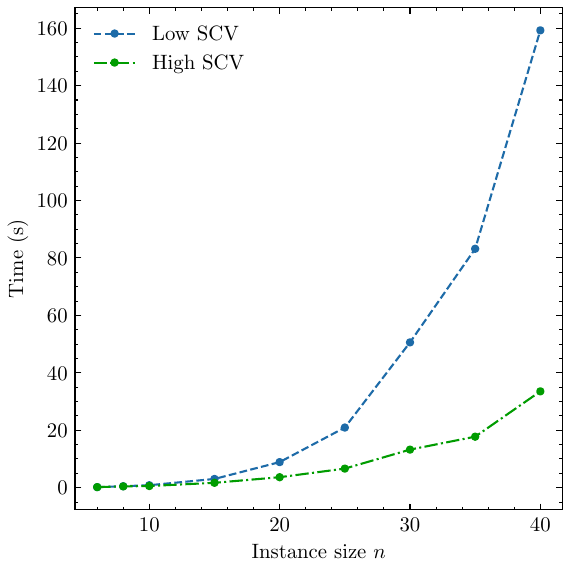}%
    \label{fig:results-runtime-optimal-schedule}}
    \caption{
      Average runtimes in seconds (A) for one LNS iteration and (B) to compute optimal appointment times.
    }
    \label{fig:results-runtime}
\end{figure*}

Figure~\ref{fig:results-runtime} visualizes the average runtimes for both a single LNS iteration ${\mathscr L}_{\rm hyb}(\Tour)$ and computing the optimal appointment times ${\mathscr L}_{\rm ph}(\Tour)$ for some tour $\Tour$. 
The runtime values were obtained by solving the benchmarking instances in the previous sections.
Figure~\ref{fig:results-runtime-lns-iteration} shows the average runtime for a single LNS iteration.
For the low SCV case, the plot shows that a single LNS iteration can take up 30 seconds for $n=40$ clients.
In view of the time limit of $960$ seconds that we used for instances with $n=40$, this means that LNS was only able to run roughly thirty iterations to solve the instance.
On the other hand, for high SCVs, the average runtime is 5 seconds for $n=40$, meaning that it would have up to 200 iterations.
Figure~\ref{fig:results-runtime-optimal-schedule} shows the average runtime for computing the optimal appointment schedule after any algorithm has decided on the final tours.
Note that this time is not included in the LNS time limit.
For high SCVs, computing the optimal schedule takes up to 40 seconds with $n=40$ clients, whereas for low SCVs it can take up to 160 seconds.
Both figures indicate that for cases with low variance, there is a very rapid increase in computational burden which renders a local search-based approach like LNS infeasible for large instances.
For cases with higher variance, the performance is considerably better, and it can consistently outperform the considered heuristics.

\subsubsection{Sensitivity analysis}
\label{sec:results-scaling-travel-times}
In Subsections~\ref{sec:results-small-instances} and~\ref{sec:results-large-instances}, we showed that MSVF is outperformed by the other algorithms, which indicates that the benchmark instances emphasize good routing decisions over appointment scheduling decisions.
We note that the impact of routing is amplified by the fact that a relatively short tour coincides with a smaller variance in travel times, in turn leading to smaller idle and waiting times. Hence, a good routing strategy provides additional benefits for appointment scheduling decisions.

In the final set of experiments, we performed a sensitivity analysis by shifting the emphasis from routing to appointment scheduling by scaling the travel times.
To this end, we used the set of benchmark instances with $n=15$ with both low and high SCV and a fixed travel time weight $\omegaT$ of $1$.
To emphasize the relative importance of appointment scheduling, we scaled the travel times with a scaling factor, for which we used values from $0.1$ to $1$ with increments of $0.1$.
The service times were unmodified.

\begin{figure}
    \centering
    \includegraphics[width=\linewidth]{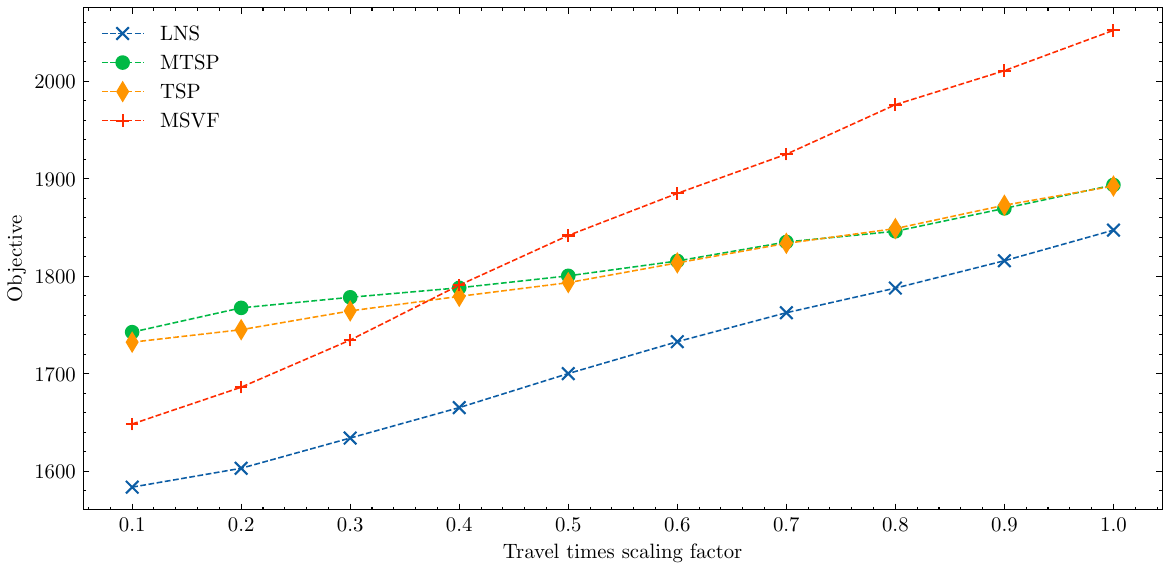}
    \caption{
    Influence of the travel time scaling factor on the average objective for each algorithm.
    }
    \label{fig:service-vs-travel}
\end{figure}

Figure~\ref{fig:service-vs-travel} shows the average objective value obtained by each algorithm as a function of the scaling factor.
As the travel-time scaling factor decreases, the relative performance of MSVF compared to the TSP-based algorithms improves, as there is more emphasis on the appointment scheduling aspect of the RAS problem.
For a scaling factor of 0.4, the MTSP, TSP, and MSVF heuristics perform about equally well; MSVF outperforms the TSP-based algorithms in case the scaling factor is smaller than 0.4.
LNS demonstrates superior performance compared to the heuristics across all instances. In particular, when the importance of routing and appointment scheduling are well balanced (with a scaling factor between 0.3 and 0.5) the gap between LNS and the best of the other algorithms is considerable. 
Overall, in situations with a weighted objective of routing and appointment scheduling, LNS is able to efficiently take both aspects into account, demonstrating the benefit of integrated decision-making.

\section{Discussion and concluding remarks}
\label{sec:conclusion}

In this paper, we have studied the integrated routing and appointment scheduling (RAS) problem, with the complicating feature that service and travel times are assumed to be of a stochastic nature. 
Our approach phrases the problem in terms of queueing-theoretic concepts, which is in contrast with the more commonly applied scenario-based approaches. 
First, we point out how the problem's objective function can be evaluated when all random variables involved are of phase-type; this is attractive because one can accurately approximate any random variables phase-type counterparts, but the downside is that the evaluation of an optimal appointment schedule can be slow, in particular for larger numbers of clients.  
This has motivated the exploration of a second class of approximations, based on the celebrated heavy-traffic approximation from queueing theory; while sacrificing some accuracy, the closed-form formulas allow extremely fast evaluation of the objective function.  
A hybrid method combines the strong aspects of both approaches: by the heavy-traffic approximation the optimal appointment times are evaluated, which are then inserted in the phase-type objective function. 
The use of this hybrid method effectively reduces the RAS problem to a more conventional routing problem, that we propose to solve using large neighborhood search. 

The effectiveness of our approach is demonstrated through an extensive set of computational experiments, with up to 40 clients.
These experiment show that it outperforms other heuristics in (almost) all instances. Overall, the gap between our approach and TSP-based heuristics tends to become smaller when the travel time is relatively large compared to the service time, the weight of the travel time is relatively large, and the variability in service times is small. 

In instances that the coefficient of variation of the service times is small and the number of clients is large (say, in the order of 40), the computation time of our approach starts to explode.
This is due to the fact that when using the hybrid method, the phase-type objective function still needs to be evaluated once per every iteration, 
which becomes prohibitively slow due to the matrix operations required (matrix exponentials, inverse matrices, etc.). 
If we wish to scale up our method so as to be able to accommodate larger numbers of clients, we need to find an alternative to the currently used phase-type evaluation of the objective function. 

We envisage great potential of our queueing-based approach in similar types of optimization problems. 
One could think of the situation of multiple service operators, where in addition to the routing and appointment decisions, one needs to determine which clients are to be served by each of the servers; one could add even one more layer in which and the number of service operators has to be chosen. 
Moreover, one could incorporate randomness in the number of clients, e.g., due to last-minute cancellations or unscheduled `high-priority' clients. 
Other possible extensions could relate to taking into account client preferences into the appointment times. 

\newpage
\appendix
\section{Phase-type distributions}\label{APP1}

Consider a discrete-time Markov chain $\{J_t\}_{t\geqslant 0}$ that is defined on the finite state space $\{1,\ldots,d+1\}$. Here, all states $1,\ldots,d$ are transient and state $d+1$ is absorbing. Then for $\{J_t\}_{t\geqslant 0}$ having an initial distribution ${\bs\alpha}\in{\mathbb R}^d$ (which is to be understood as a {\it row} vector), a phase-type random variable is defined as the time it takes to reach state $d+1.$ We then define a transition rate matrix $V\equiv (v_{kl})_{k,l=1}^d$ as follows. Denote by $v_{kl}$, with $k\not=l$ and $k,l=1,\ldots,d$, the transition rate from transient state $k$ to transient state $l$; the rate from state $k$ to the absorbing state is denoted by $v_k$, for $k=1,\ldots,d$. The diagonal elements of the matrix $V$ are chosen such that ${\bs v}=-V{\bs 1}$, where ${\bs 1}$ is a vector of all ones of appropriate dimension. A phase-type distribution is characterized through the triple $(d,{\bs\alpha}, V)$.

For such phase-type distributions, a large set of results are known; see, for instance, \cite[Section III.4]{asmussen2003}. In the context of the present paper, the most relevant result concerns the corresponding probability density function: it is given by $f(y)={\bs\alpha}\,e^{Vy}\,{\bs v}.$

Two special phase-type distributions play a crucial role in this paper: the mixture of Erlang distributions $E_K(\mu,p)$ and the hyperexponential distribution $H_2(\mu_1,\mu_2,p)$.
\begin{itemize}
    \item[$\circ$] For the mixture of Erlang distributions $E_K(\mu,p)$ we have $d=K$, ${\bs\alpha}=(1-p,p,0,\ldots,0)$ and
    \[V=\left(\begin{array}{ccccc}
    -\mu&\mu&0&\cdots&0\\
    0&-\mu&\mu&\cdots&0\\
    0&0&-\mu&\cdots&0\\
    \vdots&\vdots&\vdots&\ddots&\vdots\\
    0&0&0&\cdots&-\mu\end{array}\right).\]
    \item[$\circ$] For the hyperexponential distribution $H_2(\mu_1,\mu_2,p)$ we have $d=2$, ${\bs\alpha}= (p,1-p)$ and
    \[V=\left(\begin{array}{cc}
    -\mu_1&0\\
    0&-\mu_2\end{array}\right).\]
\end{itemize}

\section{Phase-type method for the objective function}\label{APP2}

In this appendix we point out how the objective function ${\mathscr L}(\Tour,{\bs x)}$ can be numerically evaluated given we know the phase-type description $\bar U_j$ of the random variables $U_j$; that is, we present the derivation of Theorem~\ref{thm:obj}. The procedure uses various concepts and ideas developed in \cite{wang1997}. 
As in Subsection~\ref{sec:algo_objective}, we assume that the tour $\Tour$ is given, where $j$ denotes the index of the $j$-th visited client in the given tour $\Tour$. 

The following lemma plays a crucial role. Define by $R_j$ the {\it sojourn time} of the $j$-th client, which is in this context defined as the time elapsed between the appointment time and the moment that the service of the $(j+1)$-st can start; i.e., $R_j=W_j+U_{j+1}$.
Note, importantly, that this means that in the present setting in which the server travels between locations, the {sojourn time} of the $j$-th client is {\it not} defined as the time she spends in the system (i.e., $W_i+B_i$).
For $j=1,\ldots,n$, we define
\begin{align*}{\mathscr R}_j({\bs x}) :=&\, \omegaI \,{\mathbb E}\big(x_{j} -R_{j-1}\big)\,1\{x_{j}\geqslant R_{j-1}\} +
\omegaW_j \,{\mathbb E}\big(R_{j-1} - x_{j}\big)\,1\{x_{j}< R_{j-1}\}\\
=&\,\omegaI \,{\mathbb E}\big(x_{j} - R_{j-1}\big)^+ +
\omegaW_j  \,{\mathbb E}\big(R_{j-1} - x_{j}\big)^+.
\end{align*}

\begin{lemma} For $\Tour \in{\mathscr I},{\bs x\in {\mathbb R}_+^n}$,
\[{\mathscr L}(\Tour,{\bs x)} = \omegaT
\sum_{j=1}^{n+1} {\mathbb E}\,T_{j} + \sum_{j=1}^n {\mathscr R}_j({\bs x}).\]
\end{lemma}

{\it Proof.} Recall from \eqref{objf} that the objective function is defined by
\[{\mathscr L}(\Tour,{\bs x)} = \omegaT
\sum_{j=1}^{n+1} {\mathbb E}\,T_{j} +\omegaI
\sum_{j=1}^n {\mathbb E}\,I_{j} +  \sum_{j=1}^n \omegaW_j {\mathbb E}\,W_{j}.
\]
As we have that $R_j=W_j+U_{j+1}$, the Lindley recursion becomes
\begin{align*}
    I_j= \max\{x_{j}-R_{j-1},0\},\:\:\:W_j=\max\{R_{j-1}-x_{j},0\}.
\end{align*}
The stated follows directly.\hfill$\Box$

The following lemma helps evaluating ${\mathscr R}_j({\bs x})$ when the random variables $R_{j}$ are of phase type.

\begin{lemma}\label{lemma2}
Let the non-negative random variable $X$ be of phase type of dimension $d$, with the initial distribution represented by the row vector ${\bs \alpha}\in{\mathbb R}_+^d$ and the $d\times d$ transition rate matrix given by ${V}$. Then, for any $A,B\in {\mathbb R}$ and $x\geqslant 0$,
\[A\,{\mathbb E}(X-x)^++B\,{\mathbb E}(x-X)^+= -(A+B)\,{\bs\alpha}\,V^{-1}e^{Vx}\,{\bs 1}+B\,(x+{\bs\alpha}\,V^{-1}{\bs 1}).\]
\end{lemma}

{\it Proof.}
Start by observing that
\begin{equation} \label{eqn:objIdentity}
    A(X-x)1\{X\geqslant x\}+B(x-X)1\{X<x\}=(A+B)(X-x)1\{X\geqslant x\}+B(x-X).
\end{equation}
Denote ${\bs v}:=-V{\bs 1}$. Then, using integration by parts, performing the change of variables $z:=y-x$, and denoting by $f_X(x)$ the density of $X$ (as given in Appendix \ref{APP1}),
\begin{align*}
    {\mathbb E}(X-x)^+ &= \int_x^\infty f_X(y)\,(y-x)\,{\rm d}y
    = \int_x^\infty {\bs\alpha}\, e^{Vy}{\bs v}\,(y-x)\,{\rm d}y
    \\&= -\int_0^\infty {\bs\alpha} \,V^{-1}\,e^{V(x+z)}\,{\bs v}\,{\rm d}z=\int_0^\infty {\bs\alpha} \,e^{V(x+z)}\,{\bs 1}\,{\rm d}z=-{\bs\alpha}\,V^{-1}e^{Vx}\,{\bs 1},
\end{align*}
where it is used that a matrix and its matrix exponent commute. Now using the standard result that ${\mathbb E}\,X= -{\bs\alpha}\,V^{-1}{\bs 1}$ (which of course also follows from inserting $x=0$ in the obtained expression for ${\mathbb E}(X-x)^+$), the stated follows.
\hfill$\Box$

From the above we conclude that we are left with finding the phase-type description of the random variables $R_j$. Let the phase type description $\bar U_j$ of $U_j$ be given by the the row vector ${\bs \alpha}_j\in{\mathbb R}_+^{d_j}$ and the $d_j\times d_j$ transition rate matrix given by ${V}_j$, where ${\bs v}_j:=-V_j{\bs 1}$ (where the dimension of this all-ones vector ${\bs 1}$ is $d_j$). Define the matrices $V^{(j)}$ in the following recursive manner. Let $D_j:=d_1+\ldots+d_j$ be the sum of the dimensions of phase-type random variables $\bar U_1$ up to $\bar U_j$, and let ${\mathbb O}_{m,n}$ be an all-zeroes matrix of dimension $m\times n$. With $V^{(1)}=V_1$, we recursively define
\[V^{(j)} :=\left(\begin{array}{cc}V^{(j-1)}&\tilde V_j\\{\mathbb O}_{d_j,D_{j-1}}&V_j
\end{array}\right),\]
where, for $j=2,3,\ldots$,
\[\tilde V_j=\left(\begin{array}{c} {\mathbb O}_{D_{j-2},d_j}\\-V_{j-1}{\bs 1}\,{\bs\alpha}_j\end{array}\right);\]
observe that the dimension of $-V_{j-1}{\bs 1}\,{\bs\alpha}_j$ is $d_{j-1}\times d_j$, so that
$\tilde V_j$ has dimension $D_{j-1}\times d_j$ and $V^{(j)}$ has dimension $D_j\times D_j$.

We compute the objects $F_j(y)$ and ${\bs P}_j(y)$ in the following recursive manner. Define, for $y\geqslant 0$,
\[F_j(y):= 1-{\bs P}_j(y)\,{\bs 1},\]
where ${\bs P}_1(y) :={\bs\alpha}_1e^{V^{(1)}y}$ and, for $j=2,3,\ldots$,
\[{\bs P}_j(y):= \big({\bs P}_{j-1}(x_{j}),{\bs\alpha}_j\, F_{j-1}(x_{j})\big)\, e^{V^{(j)} y},\] being a $D_j$-dimensional row vector. Mimicking the procedure developed in \cite{wang1997}, the following result can be derived.

\begin{lemma} For all $y\geqslant 0$,
\[{\mathbb P}(R_j\leqslant y) = F_j(y)= 1-{\bs P}_j(y)\,{\bs 1}.\]
\end{lemma}

The lemma directly entails that $R_j$ has a phase-type distribution of dimension $D_j$, with the initial distribution represented by the row vector 
\[{\bs\alpha}_{(j)}:=\big({\bs P}_{j-1}(x_{j}),{\bs\alpha}_j\, F_{j-1}(x_{j})\big)\in {\mathbb R}^{D_j}\] 
and the $D_j\times D_j$ transition rate matrix given by $V^{(j)}$.
Combining the above three lemmas, we find the closed form expression for the objective function ${\mathscr L}(\Tour,{\bs x)}$ as presented in Theorem~\ref{thm:obj}.
Here, we use that the component corresponding to the idle time is 
    \[ \sum_{j=1}^n  \omegaI\left(x_{j}+{\bs\alpha}_{(j)}(V^{(j)})^{-1} {\bs 1} - {\bs\alpha}_{(j)}(V^{(j)})^{-1} e^{V^{(j)}x_{j}}{\bs 1}\right), \]
whereas the component due to waiting times becomes
    \[ -\sum_{j=1}^n\omegaW_j{\bs\alpha}_{(j)}(V^{(j)})^{-1} e^{V^{(j)}x_{j}}{\bs 1}.\]

With Theorem~\ref{thm:obj} at our disposal, the numerical evaluation of ${\mathscr L}(\Tour,{\bs x)}$ has become a standard task. It requires matrix multiplication, matrix inversion, and the evaluation of matrix exponentials, which are all standard procedures in virtually any numerical software package.
For the evaluation of the matrix exponential, it is a computational advantage that the matrices $V^{(j)}$ are triangular so that the eigenvalues appear on their diagonals.

\section{Modified traveling salesman problem}

The MTSP heuristic is from \cite{zhan2021}. The idea is to modify the cost of traveling between any pair of clients by taking both the traveling cost as well as the appointment cost (idle and waiting time) into account. As a proxy, the appointment cost for some client $j$ is based on the visit of only its preceding client $i$:
\begin{equation} \label{eqn:newsvendor}
C_{ij}(x_j) = \omegaW_j {\mathbb E}(U_{ij} -x_{j})^+ + \omegaI {\mathbb E}(x_j -U_{ij})^+ ,
\end{equation}
where $U_{ij}=B_{i} +T_{ij}$.
Here $x_j$ represents the inter-appointment time between clients $i$ and $j$.
Due to the newsvendor model, we directly obtain that (\ref{eqn:newsvendor}) is minimized by \[x_j^\star=F_{ij}^{-1}\left(\frac{\omegaW_j }{\omegaW_j +\omegaI}\right),\] with $F_{ij}^{-1}(\cdot)$ denoting the inverse of $U_{ij}$. The traveling costs are then modified as
$\hat{c}_{ij} = {\mathbb E}T_{ij} + C_{ij}(x_j^\star)/\omegaT$.
First, the route ${\Tour}^{\rm MTSP}$ is determined as the solution of the asymmetric TSP with travel costs $\hat{c}_{ij}$; see Algorithm~1 of \cite{zhan2021} for an algorithmic outline. Second, the optimal appointment schedule is then obtained by minimizing the inter-appointment times
$\min_{\bs{x} \in \mathscr{X}} \mathscr{L}(\Tour^{\rm MTSP}, \bs{x})$.

It remains to specify the appointment costs $C_{ij}(x_j)$ in (\ref{eqn:newsvendor}). Due to~\eqref{eqn:objIdentity}, we may rewrite (\ref{eqn:newsvendor}) as
\[  C_{ij}(x_j) = (\omegaW_j + \omegaI)\,{\mathbb E}(U_{ij} -x_{j})^+ + \omegaI x_j -\omegaI\, {\mathbb E} U_{ij}, 
\]
such that we only need ${\mathbb E}(U_{ij} -x_{j})^+$.
As in Subsection~\ref{sec:algo_objective}, we use the phase-type description $\bar U_{ij}$ of $U_{ij}$. 
Now, for $\mbox{\scv}(U_{ij}) \ge 1$, such that $\bar U_{ij}\sim H_2(\mu_1,\mu_2,p)$ for some parameters $\mu_1,\mu_2,$ and $p$, 
\[{\mathbb E}(\bar U_{ij} -x_{j})^+ = \frac{p}{\mu_1} e^{-\mu_1 x_j} + \frac{1-p}{\mu_2}  e^{-\mu_2 x_j}. \]
Otherwise, for $\mbox{\scv}(U_{ij}) < 1$, such that $\bar U_{ij} \sim E_K(\mu,p)$ for some parameters $\mu$ and $p$, it follows after some lengthy calculations that
\[ {\mathbb E}(\bar U_{ij} -x_{j})^+= \frac{k-p-\mu x_j}{\mu (k-2)!}\Gamma(k-1, \mu x_j) + \frac{k-p}{\mu(k-1)!}(\mu x_j)^{k-1}e^{-\mu x_j}, \]
where $\Gamma(k,t)$ denotes the incomplete Gamma integral $\int_t^\infty z^{k-1} e^{-z}\,{\rm d}z$.

\section{Accuracy of hybrid approximation}\label{accht}

In this section, we consider again the appointment scheduling component of the RAS, assuming a fixed route $\Tour$. When searching for good appointment schedules, we need to find a balance between speed and accuracy, where speed is particularly important due to the iterative nature of the search procedure for the routing component.
Below, we provide a rationale for opting for the hybrid approximation over the phase-type method and the heavy traffic approximations. 

As indicated in Subsection~\ref{sec:algo_objective}, determining optimized inter-appointment times using the phase-type method ${\mathscr L}_{\rm ph}(\Tour)$ becomes prohibitively slow. 
Also, observe that ${\mathscr L}_{\rm ph}(\Tour)$ corresponds to the true optimal objective function ${\mathscr L}(\Tour) := \min_{{\bs x\in {\mathbb R}_+^n}} {\mathscr L}(\Tour,{\bs x})$ in case $U_j$ follows a mixed Erlang $E_K$ or hyperexponential $H_2$ distribution. 
Hence, we wish to investigate how the hybrid ${\mathscr L}_{\rm hyb}(\Tour)$ and heavy-traffic ${\mathscr L}_{\rm ht}(\Tour)$ approximations perform compared to the true optimum ${\mathscr L}(\Tour)$ for a variety of routes $\Tour$ in case of phase-type $U_j$.

The performance of the hybrid and heavy-traffic approximations are quantified, respectively, through the relative errors
\[
\Delta_{\rm hyb}=\frac{{\mathscr L}_{\rm hyb}(\Tour)-{\mathscr L}(\Tour)}{{\mathscr L}(\Tour)}\cdot 100\%, \:\:\:\Delta_{\rm ht} := \frac{{\mathscr L}_{\rm ht}(\Tour)-{\mathscr L}(\Tour)}{{\mathscr L}(\Tour)}\cdot 100\%.
\]
To consider a variety of instances, we determine the errors for 200 randomly generated scenarios. The mean of each of the $U_j$ is sampled uniformly from $[0.5, 1.5]$ and the {\scv} is sampled uniformly from $[0.2, 1.8]$. All instances have $40$ clients, whereas $\omegaI=0.8$ and $\omegaW_j=0.2$ for every $j$.

\begin{figure}[ht]
    \centering
    \includegraphics[width=0.57\linewidth, height=7.6cm]{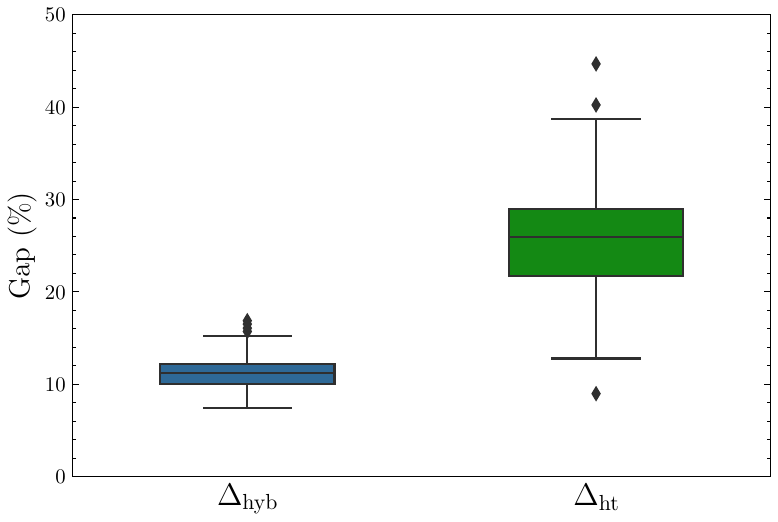}
    \caption{Box plot of $\Delta_{\rm hyb}$ and $\Delta_{\rm ht}$ for 200 random scenarios.}
    \label{fig:BP200}
\end{figure}

Figure~\ref{fig:BP200} displays the box plots of $\Delta_{\rm hyb}$ and $\Delta_{\rm ht}$. 
We see that we lose some performance, when inter-appointment times are scheduled based on heavy traffic. For instance, the median of $\Delta_{\rm hyb}$ is 11.2\%, and the median of $\Delta_{\rm ht}$ is even 25.9\%.
To assess the impact of the route $\Tour$, it is more important that the difference with ${\mathscr L}(\Tour)$ is relatively stable. After all, for the routing component, we wish to compare the absolute difference in objective function between any two alternative routes $\Tour^a$ and $\Tour^b$; this difference will not be affected by a constant error in the objective function.  
This is well achieved by ${\mathscr L}_{\rm hyb}(\Tour)$ given the narrow boxplot; the interquartile range (IQR) is only 2.18\% for $\Delta_{\rm hyb}$, whereas the IQR is 7.29\% for $\Delta_{\rm ht}$.

\begin{figure}[ht]
    \centering
    \subfloat[${\mathscr L}(\Tour)$ vs.\ ${\mathscr L}_{\rm hyb}(\Tour)$]{\label{A}{ \includegraphics[width=0.45 \linewidth]{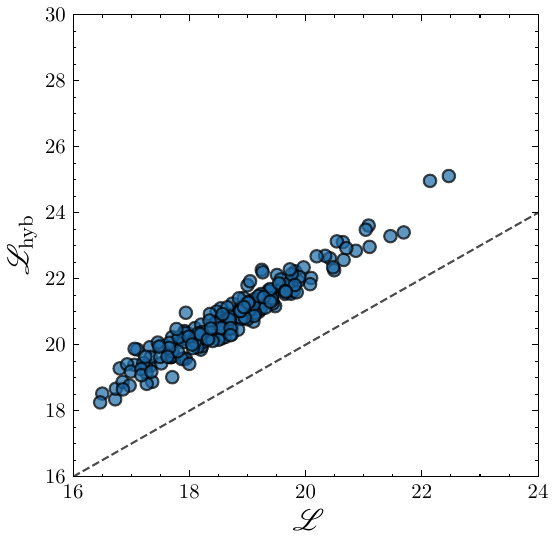} }}
    \subfloat[${\mathscr L}(\Tour)$ vs.\ ${\mathscr L}_{\rm ht}(\Tour)$]{ {\includegraphics[width=0.45\linewidth]{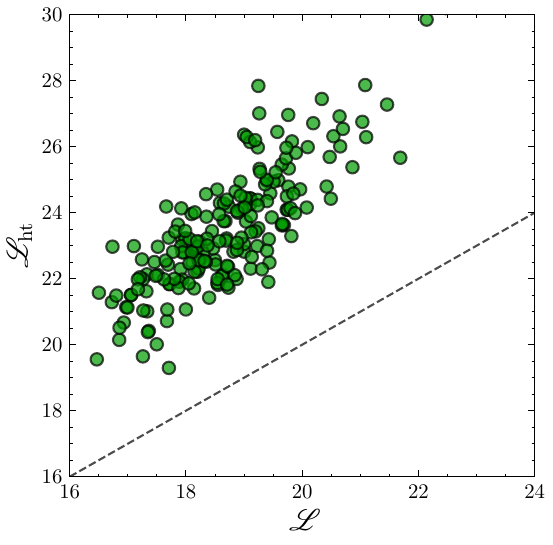} }\label{B}}
    \caption{Scatterplots of the objective functions for 200 random scenarios.}
    \label{fig: htpto}
\end{figure}


Figure~\ref{fig: htpto} presents the objective function for the 200 random scenarios using scatterplots. In line with the boxplots, Figure~\ref{A} clearly shows a strong linear relation between the hybrid approximation ${\mathscr L}_{\rm hyb}(\Tour)$ and the true optimal ${\mathscr L}(\Tour)$ with little `noise'. The relation between ${\mathscr L}_{\rm ht}(\Tour)$ and ${\mathscr L}(\Tour)$ in Figure~\ref{B} is also clearly present, but not as strong.

\begin{table}[ht]
\begin{tabular}{llccccc}
 \toprule
\multicolumn{1}{c}{No.}       & \multicolumn{1}{c}{Scenario ($\Tour$)}                     & \multicolumn{1}{c}{${\mathscr L}(\Tour)$ }           & \multicolumn{1}{c}{${\mathscr L}_{\rm hyb}(\Tour)$}             & \multicolumn{1}{c}{${\mathscr L}_{\rm ht}(\Tour)$}       & \multicolumn{1}{c}{$\Delta_{\rm hyb}$}         & \multicolumn{1}{c}{$\Delta_{\rm ht}$}   \\
\midrule
\multicolumn{1}{c}{1}         & \multicolumn{1}{l}{3 2 2 2 2 2 2 1}      & \multicolumn{1}{c}{8.35}      & \multicolumn{1}{c}{10.15}            & \multicolumn{1}{c}{11.93}      & \multicolumn{1}{c}{21.6}     & \multicolumn{1}{c}{42.8}      \\
\multicolumn{1}{c}{2}         & \multicolumn{1}{l}{1 3 2 2 2 2 2 2}      & \multicolumn{1}{c}{8.07}      & \multicolumn{1}{c}{9.54}          & \multicolumn{1}{c}{10.02}      & \multicolumn{1}{c}{18.2}  & \multicolumn{1}{c}{24.1}      \\
\multicolumn{1}{c}{3}         & \multicolumn{1}{l}{1 2 2 3 2 2 2 2}      & \multicolumn{1}{c}{8.00}      & \multicolumn{1}{c}{9.41}           & \multicolumn{1}{c}{9.59}      & \multicolumn{1}{c}{17.6}   & \multicolumn{1}{c}{21.1}      \\
\multicolumn{1}{c}{4}         & \multicolumn{1}{l}{1 2 2 2 2 3 2 2}      & \multicolumn{1}{c}{7.92}      & \multicolumn{1}{c}{9.31}           & \multicolumn{1}{c}{9.38}      & \multicolumn{1}{c}{17.6}  & \multicolumn{1}{c}{18.4}       \\
\multicolumn{1}{c}{5}         & \multicolumn{1}{l}{1 2 2 2 2 2 3 2}      & \multicolumn{1}{c}{7.83}      & \multicolumn{1}{c}{9.26}          & \multicolumn{1}{c}{9.30}      & \multicolumn{1}{c}{18.2} & \multicolumn{1}{c}{18.7}      \\
\multicolumn{1}{c}{6}         & \multicolumn{1}{l}{1 2 2 2 2 2 2 3}      & \multicolumn{1}{c}{7.65}      & \multicolumn{1}{c}{9.17}         & \multicolumn{1}{c}{9.24}      & \multicolumn{1}{c}{19.9} & \multicolumn{1}{c}{20.7}     \\
\bottomrule
\end{tabular}
\caption{\label{tab: scen} Results of different scenarios with $n= 40$ locations including depot, and their respective $\scv$'s that are arranged in batches of 5 where,  1=0.3, 2=0.9, 3=1.5, and $\omegaI=0.8$.}
\end{table}

Finally, we carried out some experiments to investigate how the objective function behaves with respect to the SVF rule. As alluded to in the introduction, the intuition behind SVF is that highly variable service times are more likely to cause (excessive) delays, which can have a `cascade' affect for all subsequent clients; scheduling clients with less variable service times earlier should result in smaller overall waiting and idle times. Again, we consider $n=40$ clients and $\omegaI=0.8$ and $\omegaW_j=0.2$ for every $j$. The 40 clients are divided into 8 batches with each batch having 5 clients with the same SCV.
In Table~\ref{tab: scen}, we considered six scenarios, where the numbers 1, 2 and 3 represent batches of clients with an SCV of 0.3, 0.9, and 1.5, respectively. Scenario 6 corresponds to SVF and has indeed the lowest objective function. In line with intuition, the objective function becomes lower when the batch of clients with the larger SCV is later in the schedule; this holds for all ${\mathscr L}(\Tour)$, ${\mathscr L}_{\rm hyb}(\Tour)$, and ${\mathscr L}_{\rm ht}(\Tour)$. Moreover, the gap $\Delta_{\rm hyb}$ for the hybrid approximation is rather stable, whereas the heavy traffic approximation has a much larger gap $\Delta_{\rm ht}$ for scenario 1.

To summarize, from our experiments we observe that the hybrid approximation is a suitable approach to compare alternative routes with optimized inter-appointment times. The variability in the heavy traffic approximation seems to large to guide the search process for the RAS.

\bibliographystyle{abbrv}
\bibliography{references}

\end{document}